\documentclass[12pt]{article}
\usepackage{amsmath,amsfonts,pinlabel,amsthm,mathptmx}

\renewcommand{\bold}[1]{\medskip \noindent {\bf \boldmath #1
                        }\nopagebreak[4]}

\newtheorem{theorem}{Theorem}[section]  

\newcommand{\half}{{\mathbb H}}
\newcommand{\integers}{{\mathbb Z}}
\newcommand{\natls}{{\mathbb N}}

\newcommand{\reals}{{\mathbb R}}

\newcommand{\makefig}[3]{
	\begin{figure}[htbp]
        \refstepcounter{figure}
	\label{#2}
        \begin{center}~
		#3~\\
		\medskip
                {Figure \thefigure.  #1}
        \end{center}
	\medskip
	\end{figure}
}

\newenvironment{pf*}[1]{%
 \begin{proof}[#1]%
}{ 
 \end{proof}
}

\newcommand{\widemargins}{
\setlength{\textwidth}{6.0in}
\setlength{\oddsidemargin}{0.25in}
\setlength{\evensidemargin}{0.25in}
}

\newcommand{\bdry}{\partial}

\newcommand{\closure}{\overline}
\newcommand{\compos}{\circ}

\usepackage{amssymb}
\newcommand{\nullset}{\varnothing}

\newcommand{\st}{\; | \;}         

\newcommand{\psibar}{{\overline{\psi}}}

\newcommand{\wt}{\widetilde}

\newcommand{\zed}{\integers}

\newcommand{\interior}{\mbox{\rm int}}

\newcommand{\Mod}{\mbox{\rm Mod}}

\newcommand{\PSL}{\mbox{\rm PSL}}

\newcommand{\Teich}{\mbox{\rm Teich}}

\newcommand{\CAT}{\mbox{\rm CAT}}

\newtheorem{prop}[theorem]{Proposition}
\newtheorem{lem}[theorem]{Lemma}
\newtheorem{cor}[theorem]{Corollary}
\newtheorem{conj}[theorem]{Conjecture}
\newtheorem{defn}[theorem]{Definition}

\newtheorem*{namedtheorem}{\theoremname}
\newcommand{\theoremname}{testing}
\newenvironment{named}[1]{\renewcommand{\theoremname}{#1}\begin{namedtheorem}}{\end{namedtheorem}}

\newcommand{\calC}{{\mathcal C}}

\newcommand{\calE}{{\mathcal E}}

\newcommand{\calG}{{\mathcal G}}

\newcommand{\calL}{{\mathcal L}}
\newcommand{\calM}{{\mathcal M}}
\newcommand{\calN}{{\mathcal N}}

\newcommand{\calP}{{\mathcal P}}

\newcommand{\calR}{{\mathcal R}}
\newcommand{\calS}{{\mathcal S}}
\newcommand{\calT}{{\mathcal T}}

\newcommand{\calV}{{\mathcal V}}

\newcommand{\glm}{{\calG \calL}}
\newcommand{\ml}{{\calM \calL}}
\newcommand{\pl}{{\calP \calM \calL}}
\newcommand{\el}{{\calE \calL}}

\usepackage{euscript}

\newcommand{\eS}{{\EuScript S}}

\widemargins

\begin{document}

\title{\bf Asymptotics of Weil-Petersson geodesics I: ending
  laminations, recurrence, and flows}
\author{\sc Jeffrey Brock, Howard Masur, and Yair Minsky\thanks{Each
    author partially 
supported by the NSF.  Jeffrey Brock was partially supported by a John Simon
Guggenheim Foundation Fellowship.}}

\maketitle

\begin{abstract}
  We define an ending lamination for a Weil-Petersson geodesic ray.
  Despite the lack of a natural visual boundary for the Weil-Petersson
  metric \cite{Brock:nc}, these ending laminations provide an
  effective boundary theory that encodes much of its asymptotic
  $\CAT(0)$ geometry.  In particular, we prove an {\em ending
    lamination theorem} (Theorem~\ref{theorem:relc}) for the
  full-measure set of rays that recur to the thick part, and we show
  that the association of an ending lamination embeds asymptote
  classes of recurrent rays into the Gromov-boundary of the curve
  complex $\calC(S)$.  As an application, we establish fundamentals of
  the topological dynamics of the Weil-Petersson geodesic flow,
  showing density of closed orbits and topological transitivity.
\end{abstract}
\tableofcontents

\newcommand{\ray}{\mathbf{r}}
\newcommand{\geod}{\mathbf{g}}
\newcommand{\elam}{\lambda}

\section{Introduction}
This paper is the first in a series considering the asymptotics of
geodesics in the Weil-Petersson metric on the Teichm\"uller space
$\Teich(S)$ of a compact surface $S$ with negative Euler
characteristic.

In many settings, measured laminations  and foliations encode the
asymptotic geometry of Teichm\"uller space.  As key examples, one has:
\begin{enumerate}
\item Thurston's natural compactification by projective measured
laminations \cite{Thurston:mcg,Bonahon:currents},
\item invariant projective measured foliations for Teichm\"uller rays
\cite{Kerckhoff:asymptotic}, and
\item the parametrization of  Bers's compactification 
by {\em end-invariants} (see \cite{Minsky:KGCC,Brock:Canary:Minsky:elc}).
\end{enumerate}
In a similar spirit, our goal will be to describe the asymptotics of
Weil-Petersson geodesics in Teichm\"uller space by the use of
laminations.  We define a notion of an {\em ending lamination} for a
Weil-Petersson geodesic ray (a geodesic from a point that leaves every
compact subset of Teichm\"uller space) and investigate its role as an
invariant for the ray.  Since the Weil-Petersson metric is not
complete, there are rays of finite
Weil-Petersson length.  These correspond to points in the
Weil-Petersson completion, which is parametrized by products of lower
dimensional Teichm\"uller spaces.  Ending laminations for such
rays are multi-curves with length functions tending to zero along
the ray.  Their initial tangents at a basepoint are dense in the unit tangent
space \cite{Brock:nc}, suggesting 
their associated multi-curves may play the role of ``rational points''
in encoding ending laminations for infinite rays.

We establish that the ending lamination is a complete asymptotic invariant for
{\em recurrent} rays, namely, those rays whose projections to the
moduli space $\calM(S) = \Teich(S)/\Mod(S)$ (the quotient of
Teichm\"uller space by the mapping class group) visit a fixed
compact set at a divergent sequence of times.  
In particular, it follows
that any two such rays starting at the same basepoint with the same
ending lamination are identical up to parametrization.  Despite the
lack of naturality described in \cite{Brock:nc}, this invariant allows
us to establish fundamentals of the topological dynamics of the
Weil-Petersson geodesic flow on the quotient $\calM^1(S) = T^1
\Teich(S) /\Mod(S)$ of the unit tangent bundle $T^1\Teich(S)$.
 We show
\begin{enumerate}
\item[(I.)] the set of closed Weil-Petersson geodesics is dense in $\calM^1(S)$
(Theorem~\ref{theorem:periodic}), and
\item[(II.)] there is a single Weil-Petersson geodesic that is dense in $\calM^1(S)$
(Theorem~\ref{theorem:dense}).
\end{enumerate} 
To the extent the ending lamination determines the ray, one can employ
properties of laminations to understand Weil-Petersson geometry.  We prove 
\newcommand{\relc}{
Let $\ray$ be a recurrent Weil-Petersson geodesic ray in
$\Teich(S)$ with ending lamination $\elam(\ray)$.  If $\ray'$
is any other geodesic ray with ending lamination
$\elam(\ray') = \elam(\ray)$ then $\ray'$ is strongly asymptotic to $\ray$.
}
\begin{theorem}{\sc (Recurrent Ending Lamination Theorem)}
\label{theorem:relc}
\relc
\end{theorem}
Here, we say $\ray$ and $\ray'$ are {\em strongly asymptotic} if there
are parametrizations for which the distance between the rays satisfies
$$\lim_{t \to \infty} d(\ray(t),\ray'(t)) = 0.$$ In particular, the
negative curvature of the Weil-Petersson metric guarantees that if
$\ray(0) = \ray'(0)$, then the rays are
identical if parametrized by arclength.

The {\em ending lamination} $\elam(\ray)$ for a ray $\ray$ arises out of
the asymptotics of simple closed curves with an a priori length bound.
Recall that by a theorem of Bers, there is a constant $L_S$ depending
only on $S$ so that for each $X \in \Teich(S)$ there is a pants
decomposition by geodesics on $X$ so that each such geodesic has
length at most $L_S$.  We call
such a $\gamma$ a {\em Bers curve} for $X$.

Given a Weil-Petersson geodesic ray $\ray$, the {\em ending
  lamination} $\elam(\ray)$ is a union of limits of Bers curves for
surfaces $X_n = \ray(t_n)$ along the ray.  In
section~\ref{section:preliminaries}, we give a precise description and
the proof that $\elam(\ray)$ is well defined.

In Proposition~\ref{prop:mu:filling}, we show that for a recurrent ray
$\ray$, the ending lamination $\elam(\ray)$ {\em fills} $S$.  Thus,
$\elam(\ray)$ determines a point in $\el(S)$, the Gromov boundary for
the curve complex $\calC(S)$ (see
\cite{Masur:Minsky:CCI,Klarreich:boundary,Hamenstaedt:boundary}).  We
remark that Theorem~\ref{theorem:relc} determines a preferred subset
\newcommand{\rel}{\calR\calE\calL} $\rel(S)\subset \el(S)$
corresponding to ending laminations for recurrent rays in the
Weil-Petersson metric.  In particular, the ending lamination
determines whether or not a ray is recurrent.

\bold{Teichm\"uller geodesics.}  The {\em Teichm\"uller metric} is a
$\Mod(S)$-invariant Finsler metric on $\Teich(S)$ measuring the
minimal 
quasi-conformal distortion of the extremal quasi-conformal mapping
between marked Riemann surfaces.  We emphasize the distinctions of our
settings from the more thoroughly studied behavior of Teichm\"uller
geodesics.

In \cite{Masur:asymptotic}, the second author shows that Teichm\"uller
geodesics with the same vertical foliation are strongly asymptotic
when the foliation is {\em uniquely ergodic} (meaning it admits a
unique transverse invariant measure), and that if a Teichm\"uller
geodesic ray is recurrent, then the vertical foliation is uniquely
ergodic.  By contrast, we note that there is no assumption of unique ergodicity for
$\elam(\ray)$ in Theorem~\ref{theorem:relc} and that 
\cite{Brock:flats} presents examples of recurrent rays with
non-uniquely ergodic laminations.
Furthermore, these examples are sharp in the sense that without the assumption of
recurrence examples are known of distinct infinite rays with the same
filling ending lamination
(see \cite{Brock:flats}, and compare \cite[\S 6]{Brock:nc}).

\bold{Visual boundaries.}  The negative curvature of the
Weil-Petersson metric (see \cite{Tromba:sectional,Wolpert:sectional})
provides for a compactification of $\Teich(S)$ by geodesic rays
emanating from a fixed basepoint $X$, the {\em visual sphere} at $X$.
Work of the first author (see \cite{Brock:nc}) demonstrates that the
compactification of $\Teich(S)$ is basepoint dependent and, moreover,
that the mapping class group fails to extend continuously to the
compactification.  

Standard arguments for topological transitivity and the density of
closed orbits that arise in Riemannian manifolds of negative curvature
involve the use of the boundary at infinity for the universal cover
and the natural extension of the action of the fundamental group to
the boundary.

The principal source of difficulty with carrying out such a line of
argument here is precisely the source of the basepoint dependence
shown in \cite{Brock:nc}.  The lack of completeness of the metric gives rise to finite-length
geodesic rays that leave every compact subset of Teichm\"uller space,
and these {\em finite rays} determine a subset of the boundary on
which the change of basepoint map is discontinuous.  While such finite
rays prevent the Weil-Petersson metric from exhibiting the more
standard boundary structure arising in the setting of Hadamard
manifolds (see \cite{Eberlein:flow}) we show the infinite length
Weil-Petersson geodesic rays determine a natural {\em boundary at
  infinity} for the Weil-Petersson completion.  
\newcommand{\vis}{ Let   $X \in \Teich(S)$ be a basepoint.
\begin{enumerate}
\item For any $Y \in \Teich(S)$ with $Y \not= X$, and any infinite ray
  $\ray$ based at $X$ there is a unique infinite ray $\ray'$ based at
  $Y$ with $\ray'(t) \in \Teich(S)$ for each $t$ so that $\ray'$ lies
  in the same asymptote class as $\ray$.
\item The change of basepoint map restricts to a homeomorphism on the infinite
 rays.
\end{enumerate}}

\begin{theorem}{\sc (Boundary at Infinity)} 
\label{theorem:vis}
\vis
\end{theorem}
Though the Weil-Petersson completion $\closure{\Teich(S)}$ of
$\Teich(S)$ does not satisfy the extendability of geodesics
requirement for a standard notion of a $\CAT(0)$ boundary to be well
defined, one can simply restrict attention to the infinite rays and
consider asymptote classes of infinite rays in the completion of the
Weil-Petersson metric, where two infinite rays are in the same {\em
  asymptote class} if they lie within some bounded Hausdorff distance
of one another.  Theorem~\ref{theorem:vis} gives a
basepoint-independent topology on these asymptote classes, and we
denote the resulting space by $\bdry_\infty \closure{\Teich(S)}$.

Any flat subspace in a $\CAT(0)$ space provides an obstruction to the
visibility property exhibited in strict negative curvature, namely,
the existence of a single bi-infinite geodesic asymptotic to any two
distinct points at infinity.  The encoding guaranteed by
Theorem~\ref{theorem:relc} of recurrent rays via laminations remedies
this conclusion to some degree, as it guarantees such a 
visibility property almost everywhere with respect to Riemannian volume measure
on the unit tangent bundle.  


\newcommand{\rvis}{
Let $\ray^+$ and $\ray^-$  be  two distinct infinite rays based at $X$.  
\begin{enumerate}
\item If $\ray^+$
is recurrent, then 
there is a single bi-infinite geodesic $\geod(t)$ so that $\geod^+ = \geod
\vert_{[0,\infty)}$ is {\em strongly} asymptotic to $\ray^+$ and 
$\geod^- = \geod \vert_{(-\infty, 0]}$ is asymptotic to $\ray^-$.
In particular, if both $\ray^+$ and $\ray^-$ are recurrent, then $\geod$ is
strongly asymptotic to both $\ray^-$ and $\ray^+$.
\item If $\mu$ in the measured lamination space $\ml(S)$ has bounded length on $\ray^\pm$ then it has
bounded length on $\geod^\pm$.
\end{enumerate}}

\begin{theorem}{\sc (Recurrent Visibility)}
\label{theorem:rvis}
\rvis  
\end{theorem}

Theorem~\ref{theorem:vis} leads one naturally to the question of
whether, as in other compactifications of Teichm\"uller space, the
laminations associated to rays serve as parameters.  Applying
Theorem~\ref{theorem:relc}, we find that such a parametrization holds
for the recurrent locus.  
\newcommand{\parametrize}{ The map $\elam$
  that associates to an equivalence class of recurrent rays its ending
  lamination is a homeomorphism to the subset $\rel(S)$ in $\el(S)$.
}
\begin{cor}
\label{cor:parametrize}
\parametrize
\end{cor}
We note that examples of \cite{Brock:flats} show this
parametrization fails in general, even when the ending lamination is
filling.

\smallskip

To describe our strategy further, we review geometric aspects of
the Weil-Petersson metric and its completion.

\bold{Weil-Petersson geometry.}  The Weil-Petersson metric $g_{\rm
  WP}$ on $\Teich(S)$ arises from the $L^2$ inner product $$\langle \varphi,\psi \rangle_{\rm WP} =
\int_X \frac{\varphi \psibar}{\rho^2} $$ on the
cotangent space $Q(X) = T_X^* \Teich(S)$ to Teichm\"uller space,
naturally the holomorphic quadratic differentials on $X$, where
$\rho(z) |dz|$
is the hyperbolic metric on $X$.

A fundamental distinction between the Weil-Petersson metric and other
metrics on Teichm\"uller space is its lack of completeness, due to
Wolpert and Chu \cite{Wolpert:noncompleteness,Chu:noncompleteness}.
It is nevertheless geodesically convex \cite{Wolpert:Nielsen}, and has
negative sectional curvatures
\cite{Tromba:sectional,Wolpert:sectional}.  

The failure of completeness corresponds precisely to {\em
pinching paths} in $\Teich(S)$ along which a simple closed geodesic on
$X$ is pinched to a cusp.  It is due to the second author that the completion
$\closure{\Teich(S)}$ is identified with the {\em augmented
Teichm\"uller space} and is obtained by adjoining noded Riemann
surfaces as limits of such pinching paths \cite{Masur:WP}.  Via this
identification, then, the
completion $\closure{\Teich(S)}$ (with its extended metric) descends to
a metric on the Mumford-Deligne compactification $\closure{\calM(S)}$
of the moduli space (cf. \cite{Abikoff:degenerating,Bers:nodes}).

The {\em Weil-Petersson geodesic flow} on the unit tangent bundle $T^1
\Teich(S)$ is the usual geodesic flow in the sense of Riemannian
manifolds with respect to the Weil-Petersson metric.  It commutes with
the isometric action of the modular group $\Mod(S)$ and defines a flow
on $\calM^1(S)$.

Because of failure of completeness, however, the geodesic flow is not
everywhere defined on $\calM^1(S)$; some directions meet the
compactification within finite Weil-Petersson distance.
The situation is remedied by the following.
\begin{prop}
\label{prop:defined}
The geodesic flow is defined for all time on a full measure subset of
$\calM^1(S)$, with respect to Liouville measure. 
\end{prop}

As a consequence, we address the question of the topological dynamics
of the geodesic flow on  $\calM^1(S)$.

The fact that 
the recurrent rays have full measure in the visual sphere allows
us to approximate directions in the unit tangent bundle arbitrarily
well by recurrent directions.  As a consequence, we have
\begin{theorem}{ \sc (Closed Orbits Dense)}
The set of closed Weil-Petersson geodesics is dense in $\calM^1(S)$.
\label{theorem:periodic}
\end{theorem}
Applying our parametrization by ending laminations of the boundary at
infinity, we may use the stable and unstable laminations for the axes
of pseudo-Anosov isometries of $\Teich(S)$ to find based at any X a
geodesic ray whose projection to $\calM(S)$ has a dense trajectory in
$\calM^1(S)$.
\begin{theorem}{\sc (Dense Geodesic)}
  There is a dense Weil-Petersson geodesic in $\calM^1(S)$.
\label{theorem:dense}
\end{theorem}

\bold{Combinatorics of Weil-Petersson geodesics.}  While the this
paper's focus on recurrence establishes the importance of the ending
lamination as a tool to analyze Weil-Petersson geodesics, it does not
directly address the connection between the combinatorics of the
lamination (in the sense of \cite{Masur:Minsky:CCII}) and the geometry
of geodesics.

We take up this discussion in
 \cite{Brock:Masur:Minsky:asymptoticsII} to prove a {\em bounded
   geometry theorem} relating bounded geometry (a lower bound for the
 injectivity radius of surfaces along the geodesic) to a bounded
 combinatorics condition analogous to bounded continued fractions, and
 vice versa. 
  These results give good control over the subset of geodesics with
  bounded geometry, and imply further dynamical consequences involving the
  topological entropy of the geodesic flow on compact invariant
  subsets.  The analogous discussion for the Teichm\"uller flow has
  been carried out by K. Rafi \cite{Rafi:Teich:short}, who obtains a
  complete description of the list of short curves along a
  Teichm\"uller geodesic in terms of the vertical and horizontal foliations.

  We expect in general that the
 ending lamination should predict extensive information
about bounded and short curves along the ray, in line with the ending
lamination theorem of \cite{Brock:Canary:Minsky:elc}.  In particular,
we make the following conjecture.
\begin{conj}{\sc (Short Curves)}
\label{conjecture:wp:comb}  
Let $\geod$ be a bi-infinite Weil-Petersson geodesic with ending
laminations $\elam^-$ and $\elam^+$ that fill the surface $S$, and let
$M \cong S \times \reals$ be a totally degenerate hyperbolic
3-manifold with ending laminations $\elam^-$ and $\elam^+$.  
Then we have
\begin{enumerate}
\item for each $\epsilon >0$ there is a $\delta>0$ so that for each
  simple closed curve
  $\gamma$ on $S$, if $\inf_t
  \ell_\gamma({\geod(t)}) < \delta$ then $\ell_\gamma(M) < \epsilon$.
\item for each $\delta' >0$ there is an $\epsilon' >0$ so that for
  each simple closed curve $\gamma$ on $S$, if $\ell_\gamma(M) <
  \epsilon'$ then $\inf_t \ell_\gamma({\geod(t)}) < \delta'$.
\end{enumerate}
\end{conj}  
Here, $\ell_\gamma(M)$ denotes the arclength of the unique geodesic
representative of $\gamma$ in $M$.
Though the present paper  will not treat them in more detail, 
hyperbolic 3-manifolds and Kleinian groups are discussed in in
\cite{Thurston:book:GTTM,Bonahon:tame,McMullen:book:RTM,Minsky:KGCC}
among other places.
The conjecture is essentially a combinatorial one, as the geometry of
hyperbolic 3-manifolds was shown to be controlled by the combinatorics
of the curve complex in
\cite{Minsky:CKGI,Brock:Canary:Minsky:elc}.

Such expected connections with ends of hyperbolic 3-manifolds motivate
other questions about the structure of the ending lamination
$\elam(\ray)$ for a ray $\ray$.
\begin{conj}
Let $\ray$ be a Weil-Petersson geodesic ray along which no simple
closed curve has length asymptotic to zero.
Then the ending lamination
$\elam(\ray)$ fills the surface.
\end{conj}
We establish this conjecture for recurrent rays in Proposition~\ref{prop:mu:filling}.

\bold{Plan of the paper.}  In section~\ref{section:preliminaries} we
set out necessary background, and give the definition of ending
lamination for a Weil-Petersson geodesic ray, establishing its basic
properties.  Section~\ref{section:flows} establishes that the geodesic
flow is defined for all time on a full measure set and gives the
natural application of the Poincar\'e recurrence theorem in this
setting.  Section~\ref{section:elc} establishes the main theorem, that
the ending lamination is a complete invariant for a recurrent ray, as
well deriving important topological properties of the ending
lamination itself that mirror the behavior of ending laminations for
hyperbolic 3-manifolds.  Finally, in section~\ref{section:dynamics} we
present applications of this boundary theory to the topological
dynamics of the Weil-Petersson geodesic flow.

\bold{Acknowledgements.}  The authors thank the Mathematical Sciences
Research Institute for its hospitality while this work was being
completed.  As this paper was in its final stages of completion, the
authors learned of an independent proof of Theorem~\ref{theorem:dense}
for dimension one Teichm\"uller spaces due to Pollicott, Weiss and
Wolpert \cite{Pollicott:Weiss:Wolpert:flow} by an explicit
construction.  We thank the referee for many useful comments as well as
suggestions for improving the exposition.

\section{Ending laminations for Weil-Petersson rays}
\label{section:preliminaries} 

In this section we begin by reviewing some of the notions and results
necessary for our discussion, provide references for background, and
give the definition of the ending lamination, establishing its basic
properties.

\bold{Teichm\"uller space and moduli space.}  The Teichm\"uller space
of $S$, $\Teich(S)$, parametrizes the marked, complete, finite-area hyperbolic
structures on $\interior(S)$, namely, pairs $(f,X)$ where 
$$f \colon \interior(S) \to X$$ is a {\em marking homeomorphism} to a
finite-area hyperbolic
surface $X$ and $(f,X) \sim (g,Y)$ if there is an isometry $\phi
\colon X \to Y$ for which $\phi \compos f$ is isotopic to $g$.  The
{\em mapping class group} $\Mod(S)$ of orientation preserving
homeomorphisms up to isotopy acts naturally on $\Teich(S)$ by
precomposition of markings, inducing an action by isometries in the
Weil-Petersson metric.  The quotient is the {\em moduli space}
$\calM(S)$, of hyperbolic structures on $\interior(S)$ (without
marking), and the Weil-Petersson metric descends to a metric on
$\calM(S)$.

\bold{Hyperbolic geometry of surfaces.}  Let $\eS$ denote the
collection of isotopy classes of essential, non-peripheral simple closed curves on
$S$. A {\em pants decomposition} $P$ is a maximal collection of distinct
elements of $\eS$ with $i(\alpha,\beta) = 0$ for any $\alpha$ and
$\beta$ in $P$.  Here,
$i \colon \eS \times \eS \to \zed$ denotes the {\em geometric intersection
number} which counts the minimal number of intersections between
representatives of the isotopy classes $\alpha$ and $\beta$ on $S$.
Given $X \in \Teich(S)$, each $\alpha \in \eS$ has a unique geodesic
representative $\alpha^*$ on $X$.  Its arclength 
determines a {\em geodesic length function}
$$\ell_\alpha \colon \Teich(S) \to \reals_+.$$
In \cite{Wolpert:Nielsen}, Wolpert proved that along a geodesic
$\geod(t)$ the length function $\ell_\alpha({\geod(t)})$ is strictly
convex.

For all that follows it will
be important to have in place the Theorem of Bers (see
\cite{Buser:book:spectra}) that given $S$, a compact orientable
surface of negative Euler characteristic, there is a constant $L_S >0$
so that for each $X \in \Teich(S)$ there is a pants decomposition
$P_X$ determined by simple closed geodesics on $X$ so that
$$\ell_\gamma(X) < L_S$$ for each $\gamma \in P_X$.  We call the pants
decomposition $P_X$ a {\em Bers pants decomposition for $X$} and the
curves in such a pants decomposition $P_X$ {\em Bers curves for $X$}.

A {\em geodesic lamination} $\lambda$ on a hyperbolic surface $X \in
\Teich(S)$ is a closed subset of $X$ foliated by simple complete
geodesics.  Employing the natural boundary at infinity for $\wt X$, a
geodesic lamination, like a simple closed curve, has a well defined
isotopy class on $X$, and we may speak of a single geodesic lamination
$\lambda$ as an object associated to $S$ with realizations on each
hyperbolic structure $X \in \Teich(S)$ (see
\cite[Ch. 8]{Thurston:book:GTTM}, \cite{Hatcher:ml}, and \cite{Bonahon:currents}).  The
realizations of geodesic laminations on $X$ may be given the Hausdorff topology, and
the correspondence between realizations of $\lambda$ on different
surfaces $X$ and $X'$ gives a homeomorphism.  Hence, we refer to a
single {\em geodesic lamination space} $\glm(S)$.

A geodesic lamination $\lambda$ equipped with a {\em transverse
  measure} $\mu$, namely a measure on each arc transverse to the
leaves of $\lambda$ invariant under isotopy preserving intersections
with $\lambda$, determines a {\em measured lamination}.  The
lamination $\lambda$ is called the {\em support} of the measured
lamination $\mu$ and is denoted by $|\mu|$.  The simple closed curves
with positive real weights play the role of Dirac measures, and the
measured lamination space $\ml(S)$ 
is identified with the closure of the image of the embedding
$\iota \colon \eS \times \reals_+ \to \reals_+^\eS$ by 
$$\langle \iota(\alpha, t) \rangle_\beta = t \cdot i(\alpha,\beta).$$
(see
\cite{FLP,Thurston:book:GTTM,Bonahon:currents}).  The natural
action of $\reals^+$ on $\ml(S)$ by scalar multiplication of  transverse
measures gives rise to 
Thurston's projective measured lamination space $\pl(S) = (\ml(S)
- \{0\}) / \reals_+$.  Throughout, $[\mu]$ will denote the
projective class in $\pl(S)$ of a nonzero measured lamination $\mu \in
\ml(S)$.

The geodesic length function for a simple closed curve extends 
to a bi-continuous
function
$$\ell_{.}(.) \colon \ml(S) \times \Teich(S) \to \reals_+$$
by defining $\ell_{t\cdot \gamma}(X) = t (\ell_\gamma (X))$ for $t \in
\reals_+$ and $\gamma \in \eS$, and setting
$$\ell_\mu(X) = 
\lim_{n \to \infty} s_n (\ell_{\gamma_n}(X))$$ (see
\cite{Kerckhoff:eq:analytic,Bonahon:currents}).  Wolpert strengthens
his convexity result for simple closed curves to apply to this ``total
length'' of a measured  
lamination (see  \cite{Wolpert:glf})
\begin{theorem}[Wolpert]  Given a Weil-Petersson gedesic $\geod(t)$, the length $\ell_\mu({\geod(t)})$ of a
  measured lamination is a strictly convex function of $t$.
\label{theorem:strictly}
\end{theorem}

\bold{Curve and arc complexes.}  The complex of curves $\calC(S)$
associated to the surface $S$ is a simplicial complex whose vertices
are elements of $\eS$,
and whose $k$-simplices span $k+1$-tuples of vertices
whose corresponding isotopy classes can be realized as a pairwise
disjoint collection of simple closed curves on $S$.  By convention, we
obtain the {\em augmented curve complex} by adjoining the empty
simplex and denote
$$\widehat{\calC(S)} = \calC(S) \cup \nullset.$$

It was shown in \cite{Masur:Minsky:CCI} that the curve complex
$\calC(S)$ is a $\delta$-hyperbolic path metric space.  Any such space
carries a natural {\em Gromov boundary}, which is identified with
asymptote classes of quasigeodesic rays
where two rays are
{\em asymptotic} if they lie within uniformly bounded Hausdorff
distance.  
Klarreich showed \cite{Klarreich:boundary} (see also
\cite{Hamenstaedt:boundary}) that the Gromov boundary is identified
with the space $\el(S)$ of geodesic laminations that arise
as supports  of {\em filling} measured laminations.  (A lamination
$\mu\in \ml(S)$ is filling if every simple 
closed curve $\gamma$ satisfies $i(\mu,\gamma) >0$).  The space
$\el(S)$ inherits the quotient topology from $\ml(S)$, but it is a
Hausdorff subspace of this quotient; this topology is sometimes called
the {\em measure-forgetting topology} or the {\em Thurston topology}
\cite{Canary:Epstein:Green}.

\newcommand{\collar}{\mathbf{collar}}
Given a reference hyperbolic structure $X \in \Teich(S)$, for each $\gamma \in
\eS$
there is a $\delta_\gamma >0$ such that the neighborhood
$\calN_{\delta_\gamma}(\gamma^*)$   of the geodesic representative
$\gamma^*$ on  $X$ 
is a regular neighborhood, and
so that for each $\eta$ with
$i(\eta,\gamma) = 0$, we have disjoint neighborhoods
$$\calN_{\delta_\gamma} (\gamma^*) \cap \calN_{\delta_\eta}(\eta^*)=\emptyset.$$
Given a simplex $\sigma \subset \calC(S)$, denote by
  $\sigma^0$ the set of vertices of $\sigma$ and let $\collar(\sigma)$ be
the union 
   $$\bigcup_{\gamma \in \sigma^0} \calN_{\delta_\gamma}(\gamma^*).$$  
Fixing this notation, we make the following definition.
\begin{defn}
  Let $\lambda$ be a connected geodesic lamination.  The {\em
    supporting subsurface} $S(\lambda) \subset S$ is the compact
  subsurface up to isotopy represented by the smallest subsurface $Y
  \subset X$ containing the realization of $\lambda$ as a geodesic
  lamination on $X$, whose non-peripheral boundary curves are a union
  of curves in $\collar(\sigma)$ for some $\sigma \in
  \widehat{\calC(S)}$.
\end{defn}

\bold{The pants complex.}  A quasi-isometric model was obtained for
the Weil-Petersson metric in \cite{Brock:wp} using pants
decompositions of surfaces.  We say two pants decompositions $P$ and
$P'$ are related by an {\em elementary move} if $P'$ is obtained from
$P$ by replacing a curve $\alpha$ in $P$ with a curve $\beta$ in such
a way that $i(\alpha,\beta)$ is minimized.  Let $P(S)$ denote the
graph whose vertices represent distinct isotopy classes of pants
decompositions of $S$, or maximal simplices in $\calC(S)$, and whose
edges join vertices that differ by an elementary move.

Hatcher and Thurston showed that $P(S)$ is connected (see
\cite{Hatcher:pants}) so we may consider the edge metric on $P(S)$ as
a distance on the pants decompositions of $S$.  Letting $Q \colon
P(S) \to \Teich(S)$ be any 
map that associates to $P$ a surface $X$ on which $P$ is a Bers pants
decomposition. 
\begin{theorem}{\rm (\cite[Thm. 1.1]{Brock:wp})}
The map $Q$ is a quasi-isometry.
\label{theorem:pants:quasi}
\end{theorem}
In other words, the map $Q$ distorts distances by a bounded
multiplicative factor and a bounded additive constant.  

\bold{The Weil-Petersson completion and its strata.}  The non-completeness
of the Weil-Petersson metric corresponds to finite-length paths in
Teichm\"uller space along which length functions for simple closed
curves converge to zero.  In
\cite{Masur:WP}, the completion is described  concretely as the {\em
augmented Teichm\"uller space} \cite{Bers:nodes,Abikoff:degenerating}
obtained from Teichm\"uller space by adding strata consisting of
spaces $\calS_\sigma$ defined by the vanishing of length functions
$$\ell_\alpha \equiv 0 $$ for each $\alpha \in \sigma^0$ where $\sigma$ is a
simplex in the augmented curve complex $\widehat{\calC(S)}$.  Points in the $\sigma$-null
strata $\calS_\sigma$ correspond to {\em nodal} Riemann surfaces $Z$,
where (paired) cusps are introduced along the curves in $\sigma^0$.

One can describe the topology via extended Fenchel Nielsen
coordinates:
Given a pants decomposition $P$, the usual coordinates map $\Teich(S)$
to $\prod_{\gamma\in P}\reals\times\reals^+$, where the first
coordinate of each pair measures twist 
and the second is the length function of the corresponding vertex of $P$.
We extend this to allow length 0, and take the quotient by identifying 
$(t,0)\sim(t',0)$ in each $\reals\times\reals^+$ factor. The topology
near any point of a stratum $\calS_\sigma$, where $\sigma^0\subset P$,
is such that this map is a homeomorphism near that point.

Then the strata $\calS_\sigma$ are naturally products of lower
dimensional Teichm\"uller spaces corresponding to the complete, finite-area
hyperbolic ``pieces'' of the nodal surface $Z \in \calS_\sigma$.

As observed in \cite{Wolpert:compl,Masur:Wolf:WP} 
The completion $\closure{\Teich(S)}$ has the structure of a $\CAT(0)$
space: it is a length space, satisfying the sub-comparison property
for chordal distances in comparison triangles in the Euclidean plane
(see \cite[II.1, Defn. 1.1]{Bridson:Haefliger:npc}).  Given $(X,Y) \in
\closure{\Teich(S)} \times \closure{\Teich(S)}$ we will denote by
$\overline{XY}$ the unique Weil-Petersson geodesic joining $X$ to $Y$.  
Then the {\em main stratum}, $\calS_\nullset$, is simply the full
Teichm\"uller space $\Teich(S)$.

Apropos of this convention, we recall the fundamental non-refraction
for geodesics on the Weil-Petersson completion.
\begin{theorem}[\cite{Daskalopoulos:Wentworth:mcg,Wolpert:compl}]{\sc
    (Non-Refraction in the Completion)} 
  Let $\overline{XY}$ be the geodesic joining $X$ and $Y$ in
  $\closure{\Teich(S)}$, and let $\sigma_-$ and $\sigma_+$ be the
  maximal simplices in the curve complex so that $X \in
  \calS_{\sigma_-}$ and $Y \in \calS_{\sigma_+}$.  If $\eta = \sigma_-
  \cap \sigma_+$, then we have
 $$\interior(g) \subset \calS_\eta.$$
\label{theorem:non-refraction}
\end{theorem}
We remark that in the special case that $X$ and $Y$ lie in the
interior of Teichm\"uller space the theorem is simply a restatement of
Wolpert's geodesic convexity theorem (see \cite{Wolpert:Nielsen}).  A
consequence of Theorem~\ref{theorem:non-refraction} is a
classification of elements of $\Mod(S)$ in terms of their action by
isometries of the Weil-Petersson completion $\closure{\Teich(S)}$
(see \cite{Daskalopoulos:Wentworth:mcg,Wolpert:compl}).  In particular, a
mapping class $\psi$ is {\em pseudo-Anosov} if no non-zero power of
$\psi$ preserves any isotopy class of simple closed curves on $S$.  As
in the setting of the Teichm\"uller metric, $\psi$ preserves an
invariant Weil-Petersson {\em geodesic axis} $A_\psi \subset
\Teich(S)$ on which it acts by translation.

\bold{Weil-Petersson geodesic rays  and ending laminations.}  
Allowing $\omega = \infty$, a Weil-Petersson geodesic ray 
 is a geodesic
$$\ray \colon [0, \omega) \to \Teich(S)$$ 
parametrized by arclength, so that $\ray(t)$ leaves every compact
subset of Teichm\"uller space. Note this means that even when
$\omega <\infty$, the ray cannot be extended further.

Although triangles in a $\CAT(0)$ space can fail the stronger
thin-triangles condition of Gromov hyperbolicity, the comparison
property for triangles suffices to guarantee that there is still a
well defined notion of an {\em asymptote class} for a geodesic ray:
two rays $\ray$ and $\ray'$ lie in the same asymptote class, or are
{\em asymptotic} if there is a $D >0$ so that $$d(\ray(t),\ray'(t)) <
D$$ for each $t$.

Fixing a basepoint $X \in \Teich(S)$, however, it is natural in the
setting of negative curvature to consider the sphere of geodesic rays
emanating from $X$ ($\ray(0)=X$), which we denote by $\calV_X(S)$, or the Weil-Petersson
{\em visual sphere}.  Geodesic convexity (see \cite{Wolpert:Nielsen})
guarantees that we can compactify Teichm\"uller space by appending
$\calV_X(S)$. 

We call a simple closed curve $\gamma \in \eS$ a {\em Bers curve for
  the ray $\ray$} if there is a $t \in [0,\omega)$ for which $\gamma$
is a Bers curve for $\ray(t)$.

We associate a geodesic lamination $\elam(\ray)$ to a
ray $\ray$ as follows.
\begin{defn}
  An {\em ending measure} for a geodesic ray $\ray(t)$ is
  any representative $\mu\in\ml(S)$ of a limit of projective classes $[\gamma_n]\in\pl(S)$
  of any infinite sequence of distinct Bers curves for $\ray$.
\end{defn}

\bold{Remark.} The definition of ending measures parallels Thurston's
definition of the ending lamination for a simply degenerate end of a
hyperbolic 3-manifold  (see
\cite[Ch. 9]{Thurston:book:GTTM}).  We remark that is possible for an
ending measures to be supported on a subsurface of $R \subset S$,
while the geometry of the complement of $R$ stabilizes along $\ray$.
This explains the use of infinite sequences of Bers curves rather than
Bers pants decompositions in the definition, since these may intersect
non-trivially in a subsurface whose geometry is converging.

Given $L >0$ there may be a fixed curve $\gamma$ that satisfies
$\ell_\gamma({\ray(t)}) \le L$ for each $t$.  Those $\gamma$ that have no
positive lower bound to their length, however, play a special role.
\begin{defn}
  A simple closed curve $\gamma$ is a {\em pinching curve} for $\ray$ if
  $\ell_\gamma({\ray(t)}) \to 0$ as $t \to \omega$.
\end{defn}

A single ray can exhibit both types of behavior, motivating the
following definition.
\begin{defn}
  If $\ray(t)$ is a Weil-Petersson geodesic ray, the {\em ending
  lamination} $\elam(\ray)$ for $\ray$ is the union of the pinching
  curves and the geodesic laminations arising as supports of ending
  measures for $\ray$.
\end{defn}

To justify the definition we must show that pinching curves and
supports of ending measures together have the underlying structure of
a geodesic lamination.  Specifically, we must show that pinching curves
and ending measures have no transverse intersections, or that
$i(\mu_1,\mu_2) = 0$ for any pair of pinching curves or ending
measures.

We first establish the following basic property of
ending measures.
\begin{lem}
If a $\ray$ has finite length then its collection of ending measures
is empty.
\end{lem}
\begin{proof}
It suffices to show that 
if $\ray$ has finite length then 
there does not exist an infinite
sequence of distinct Bers curves.  

But a finite-length ray $\ray(t)$ converges to a nodal surface $Z$ in the Weil-Petersson
completion $\overline{\Teich(S)}$, and for each simple closed curve $\gamma$ on
$S$ either 
\begin{enumerate}
\item  there is a pinching curve $\alpha$ for which $i(\alpha,\gamma)
  > 0$, or
\item the length of $\gamma$ converges along the ray $\ray(t)$ to its
  length on $Z$.
\end{enumerate} In the first case, the length of $\gamma$ on $\ray(t)$
diverges as $t \to \omega$ by the collar lemma (see \cite{Buser:book:spectra}).
It follows that the union of Bers curves over all surfaces $\ray(t)$ is
finite.
\end{proof}

Theorem~\ref{theorem:strictly} guarantees that each pinching curve
$\gamma$ for $\ray$ has length decreasing in $t$.  By showing their
boundedness along infinite rays, we may apply
Theorem~\ref{theorem:strictly} again to see the same holds for ending
measures.

\begin{lem}
\label{lemma:ending:bounded}
Let $\mu$ be any ending measure for $\ray$.
Then $\ell_\mu({\ray(t)})$ is decreasing in $t$.
\end{lem}

\begin{proof}Assume $\ray$ is based at $X \in \Teich(S)$.  Let
  $\gamma_n$ be a sequence of Bers curves for the ray $\ray$ so that
  the length of $\gamma_n$ is infimized at $\ray(t_n)$, and for which
  $t_i < t_{i+1}$, $i \in \natls$.  Let $[\mu]$ be any accumulation
  point of the sequence of projective classes $[\gamma_n]$ in
  $\pl(S)$.  Then $\mu$ is an ending measure for $\ray$.  We may
  assume, after rescaling, that $\mu$ is the representative in the
  projective class $[\mu]$ with $\ell_\mu(X) = 1.$

Letting $s_n >0$ be taken so that $$s_n =
\frac{1}{\ell_{\gamma_n}(X)},$$ the measured laminations $s_n
\gamma_n$ satisfy $\ell_{s_n \gamma_n}(X) = 1$ for each $n$, and it
follows that $s_n \gamma_n \to \mu$ in $\ml(S)$.

Fixing a value $t' > 0$, there is an $N'$ so that for $n > N'$, we
have $t_n > t'$.  Applying strict convexity of the length function
$\ell_{\gamma_n}({\ray(t)})$ as a function of $t$,
\cite{Wolpert:Nielsen}, we conclude that 
$$\ell_{s_n \gamma_n}({\ray(t')})  < 1 $$ for each $n > N'$.  
We conclude that $$\ell_\mu({\ray(t')}) \le 1.$$ Since $t' >0$ was
arbitrary, and $\ell_\mu({\ray(t)})$ is a strictly convex function of $t$
by Theorem~\ref{theorem:strictly},
we conclude that $\ell_\mu({\ray(t)})$
is decreasing in $t$.
\end{proof}

For future reference, we establish the following continuity property
for the behavior of bounded length laminations along rays.
\begin{lem}
Let $\ray_n \to \ray$ be a convergent sequence of rays in the visual sphere
$\calV_X(S)$.  Then if $\mu_n$ is any sequence of ending measures or
weighted pinching curves for $\ray_n$, any representative $\mu \in
\ml(S)$ of the limit $[\mu]$ of projective classes $[\mu_n]$ in $\pl(S)$ has
bounded length along the ray $\ray$.
\label{lem:continuity}
\end{lem}

\begin{proof}
After normalizing so that $\ell_{\mu_n}(X) = 1$ we may assume that 
$\ell_{\mu_n}({\ray_n(t)}) \le 1$ along $\ray_n$.
Then for each surface $Y = \ray(s)$ along $\ray$
there are surfaces $X_n = \ray_n(s)$ with $X_n \to Y$ in $\Teich(S)$.  Then
$\ell_{\mu_n}({X_n}) \to \ell_\mu(Y)$
and thus we have $\ell_\mu(Y) \le 1$.  Since $s$ is arbitrary, the
Lemma follows. 
\end{proof}

\begin{prop}
  Given a ray $\ray$, the union $\elam(\ray)$ is a non-empty
  geodesic   lamination.
\label{prop:lambda:defined}
\end{prop}

\begin{proof} We first show that given $\ray$, there exists either a
  pinching curve or an ending measure for $\ray$.  If $\ray$ is a ray
  of finite length, then it terminates in the completion at a nodal
  surface $Z$ in a boundary stratum $\calS_\sigma$.  It follows that
  each curve $\gamma$ associated to a vertex of $\sigma$ has length
  tending to zero along $\ray$ and is thus a pinching curve for
  $\ray$.

  Assume there are no pinching curves for $\ray$.  Then, since
  $\ray$ leaves every compact subset of $\Teich(S)$, and it does not
  terminate in the completion, it follows that it has infinite
  Weil-Petersson length.  Then we claim there is a non-zero ending
  measure $\mu$ for $\ray$.  It suffices to show that there are
  infinitely many distinct Bers curves $\gamma_n$ for surfaces
  $\ray(t_n)$, with $t_n \to \infty$.  But otherwise, the set of all
  Bers pants decompositions along the ray is also finite.  By
  Theorem~\ref{theorem:pants:quasi}, we obtain a
  bound for the length of the ray $\ray$ via the quasi-isometry $Q$,
  contradicting the assumption that $\ray$ was infinite.

As in the definition of the ending lamination for hyperbolic
3-manifolds \cite[Ch. 8]{Thurston:book:GTTM}, it suffices to
show that for any pair  $\mu_1$ and $\mu_2$ of weighted pinching
curves or ending measures, that the intersection number satisfies
$$i(\mu_1,\mu_2) = 0.$$ We note first that by the collar lemma any two
pinching curves for $\ray$ must be disjoint.  Furthermore, if
$\gamma$ is a pinching curve for $\ray$, 
then $\gamma$ is disjoint from each Bers curve on $\ray(t)$ for
$t$ sufficiently large.  Thus, if $\mu$ is an ending measure for
$\ray(t)$, then we have $i(\gamma,\mu) = 0$ as well.
Thus we reduce to the case that $\mu_1$ and $\mu_2$ are both ending
measures.

Assume that $i(\mu_1,\mu_2) > 0$.  We note in particular that if
$\mu_1$ and $\mu_2$ fill the surface, Lemma~\ref{lemma:ending:bounded}
guarantees that 
the ray $\ray(t)$ defines a path of
surfaces that range in a compact family in $\Teich(S)$ by Thurston's
{\em Binding Confinement} (see \cite[Prop. 2.4]{Thurston:hype2}).
This contradicts the assumption that $\ray$ leaves every compact
subset of $\Teich(S)$.

More generally, let $\mu_1$ and $\mu_2$ fill a proper essential
subsurface $Y \subset S$.  Then a more general version of binding
confinement, {\em Converge on Subsurface} (see
\cite[Thm. 6.2]{Thurston:hype2}), together with
Lemma~\ref{lemma:ending:bounded} ensures that the representations
$\rho_t \colon \pi_1(S) \to \PSL_2(\reals)$ for which $\ray(t) = \half^2
/\rho_t(\pi_1(S))$ have restrictions to $\pi_1(Y)$ that converge up to
conjugacy after passing to a subsequence.

It follows, for any curve
$\eta\in\calC(Y)$, that the length $\ell_\eta({\ray(t)})$ is bounded
away from zero and infinity. 
In particular $Y$ contains no pinching curves. 
By the collar lemma $\eta$ has a collar neighborhood of definite
width in each $\ray(t)$, which implies for any sequence $\gamma_n$ of
Bers curves on $\ray(t_n)$, that $i(\eta,\gamma_n)$ is bounded above. 

Each ending measure $\mu_i$ is a limit of weighted Bers curves $s_n
\gamma_n$, with $s_n \to 0$, so it follows that for each $\eta \in
\calC(Y)$ we have 
$$i(\eta,\mu_i) = \lim_{n\to\infty} i(\eta,s_n \gamma_n) = 0,$$ which contradicts
that the support of $\mu_i$ intersects $Y$.

We conclude that $i(\mu_1,\mu_2) = 0$, and thus that the set of
complete geodesics in the support of all pinching curves and ending
measures forms a closed subset consisting of disjoint complete geodesics,
namely, a geodesic lamination.
\end{proof}

We note the following corollary of the proof.
\begin{cor}
  Let $\mu \in \ml(S)$ be any lamination whose length is bounded along
  the ray $\ray$.  Then if $\mu'$ is an ending measure
  for $\ray$ or a measure on any pinching curve for $\ray$, we have
  $$i(\mu,\mu') = 0.$$
\label{cor:bounded:length}
\end{cor}
\begin{proof}
  The proof of Proposition~\ref{prop:lambda:defined} employs only the
  bound on the length of $\mu_1$ and $\mu_2$ along the ray to show the
  vanishing of their intersection number.  The argument applies
  equally well under the assumption that $\mu_1$ is a simple closed
  curve of bounded length, and $\mu_2$ is an ending measure, or a
  weighted pinching curve.  Letting $\mu$ play the role of $\mu_1$ and
  $\mu'$ play the role of $\mu_2$, the Corollary follows.
\end{proof}

By Thurston's classification of elements of $\Mod(S)$, a pseudo-Anosov
element $\psi \in \Mod(S)$ determines 
laminations $\mu^+$ and $\mu^-$ in $\ml(S)$, invariant by $\psi$ up to
scale \cite{Thurston:mcg}.  
Each determines an unique projective class in $\pl(S)$, the so-called
{\em stable} and {\em unstable} laminations for $\psi$, and arises as
a limit of iteration of $\psi$ on $\pl(S)$.  Specifically, given a
simple closed curve $\gamma$, we have
$$[\mu^+]  = \lim \psi^n([\gamma]) \ \ \ \ \text{and} \ \ \ \
[\mu^-] = \lim \psi^{-n}([\gamma])$$ in $\pl(S)$.  Similarly, each $X
\in A_\psi$, the axis of $\psi$, determines a {\em forward ray} $\ray^+$ based at $X$ so
that $\psi(\ray^+) \subset \ray^+$ and a {\em backward ray} $\ray^-$
at $X$ so that $\ray^- \subset \psi(\ray^-).$ Invariance of the axis
$A_\psi$, then, immediately gives the following relationship between
the stable and unstable laminations for $\psi$ and the ending
laminations for the forward and backward rays at $X$ for the invariant
axis $A_\psi$.
\begin{prop}
Let $\psi \in \Mod(S)$ be a pseudo-Anosov element with invariant axis
$A_\psi$.  Let $X \in A_\psi$, and let $\ray^+$ and $\ray^-$ be
the forward and backward geodesic rays at $X$ determined by
$A_\psi$.  Then we have
$$|\mu^+| = \elam(\ray^+) \ \ \ \ \text{and} \ \ \ \
|\mu^-| = \elam(\ray^-)$$ where $\mu^+$ is the stable lamination for
$\psi$ and $\mu^-$ is the unstable lamination.
\label{prop:stable:unstable}
\end{prop}
\begin{proof}
  Letting $\gamma$ be a Bers curve for the surface $X$, the projective
  class $[\mu^+]$ of $\mu^+$ is the limit of the projective classes
  $[\gamma_n]$ where $\gamma_n = \psi^n(\gamma)$ and likewise,
  $[\mu^-]$ is the limit of $[\gamma_{-n}]$.  Since $\gamma_n$ is a
  Bers curve for $\psi^n(X)$, it follows that $\mu^+$ and $\mu^-$ are
  ending measures $\ray^+$ and $\ray^-$, respectively.  Since $\mu^+$
  fills the surface, any other ending measure $\mu$ for $\ray^+$ has
  intersection number $i(\mu^+,\mu) = 0$, so we have $\elam(\ray^+) =
  |\mu^+|$ and likewise $\elam(\ray^-) = |\mu^-|$.
\end{proof}

\section{Density, recurrence, and flows}
\label{section:flows}
This section establishes fundamentals of the Weil-Petersson geodesic
flow on $\calM^1(S)$, which, while standard for complete Riemannian
manifolds of negative curvature, require more care due to the lack of
completeness of the Weil-Petersson metric.  In particular, the
non-refraction of geodesics at the completion,
Theorem~\ref{theorem:non-refraction} plays a crucial role in
establishing that (1) the bi-infinite and recurrent geodesics each
have full measure (Proposition~\ref{prop:defined} and~\ref{prop:poincare}), and (2) each asymptote class of infinite rays has a
representative based at each $X \in \Teich(S)$ (Theorem~\ref{theorem:vis}).

\medskip

In \cite{Brock:nc}, the $\CAT(0)$ geometry of the Weil-Petersson
completion and Theorem~\ref{theorem:non-refraction} are employed to
show the following.
\begin{theorem}{\rm (\cite[Thm. 1.5]{Brock:nc})}
The finite rays are dense in the visual sphere. 
\end{theorem}

Wolpert observed that one obtains the following generalization (see
\cite[Sec. 5]{Wolpert:compl}).
\begin{theorem}[Wolpert]
Restrictions to $\Teich(S)$ of Weil-Petersson geodesics in
$\closure{\Teich(S)}$ 
joining pairs of maximally noded surfaces are dense in the unit
tangent bundle $T^1\Teich(S)$.  
\end{theorem}

We recall a key element of the proof.
\begin{lem}[Wolpert]
The finite rays have measure zero in the visual sphere.
\end{lem}
\noindent(See \cite{Wolpert:compl,Wolpert:glf}).

\begin{proof}
  Given a simplex $\sigma$ in $\calC(S)$, consider the natural geodesic
  retraction map from a given null-stratum $\calS_\sigma$ onto the
  unit tangent sphere at $X \in \Teich(S)$, sending each point $Z \in
  \calS_\sigma$ to the unit tangent at $X$ in the direction of the
  unique geodesic from $X$ to $Z$.  Wolpert observes this map is
  Lipschitz from the intrinsic metric on $\calS_\sigma$ to the
  standard metric on the unit tangent sphere.  As each stratum has
  positive complex co-dimension, the image of $\closure{\Teich(S)}
  \setminus \Teich(S)$ has Hausdorff measure zero in the (real
  co-dimension 1) visual sphere.  It follows that {\em infinite}
  directions have full measure.
\end{proof}
Proposition~\ref{prop:defined}
follows as an immediate corollary.
\begin{named}{Proposition~\ref{prop:defined}}
  The geodesic flow is defined for all time on a full Liouville
  measure subset $\calM^1_\infty(S)$ of $\calM^1(S)$, consisting of
  bi-infinite geodesics.
\end{named}
\begin{proof}
  That the infinite rays have full-measure in the unit tangent bundle
  $T^1_X \Teich(S)$ at $X \in \Teich(S)$ implies that the directions
  determining bi-infinite geodesics have full measure in
  $T^1_X\Teich(S)$.  By Fubini's theorem, the union over $X$ of their
  projections determines a flow-invariant subset of $\calM^1(S)$ of
  full measure.
\end{proof}

A geodesic ray $\ray$ based at $X \in \calM(S)$ is {\em divergent} if
for each compact set $K \subset \calM(S)$, there is a $T$ for which
$\ray(t) \cap K = \nullset$ for each $t >T$.  A ray $\ray$ is called
{\em recurrent} if it is not divergent.  

Alternatively, Mumford's compactness
theorem \cite{Mumford:compact}, guarantees that given $\epsilon > 0$
the ``$\epsilon$-thick-part''
$$\Teich_{\ge \epsilon}(S) = \{ X \in \Teich(S) \st  \ell_\gamma(X) \ge
\epsilon, \ \gamma \in \eS
\}$$ of Teichm\"uller space projects
to a compact subset of $\calM(S)$.  Thus we may characterize recurrent
rays equivalently by the condition that there is an $\epsilon >0$ and
a sequence of times $t_n \to \infty$ so that $\ray(t_n) \subset
\Teich_{\ge \epsilon}(S)$.

A geodesic $\geod$ is {\em
  doubly recurrent} if basepoint $X \in \geod$ divides $\geod$ into
two recurrent rays based at $X$.

Taking Proposition~\ref{prop:defined} together with the Poincar\'e
recurrence theorem, we have the following.
\begin{prop}
\label{prop:poincare}
  The recurrent rays and doubly recurrent geodesics in $\calM^1(S)$ 
  determine full-measure invariant subsets.
\end{prop}

\begin{proof}
  The geodesic flow is volume-preserving on $\calM^1(S)$, by Liouville's theorem (see
  \cite[\S2, Thm. 2]{Cornfeld:Fomin:Sinai:ergodic}), and thus
  finiteness of the Weil-Petersson volume of $\calM(S)$
  (\cite{Masur:WP,Wolpert:finite}), and hence of $\calM^1(S)$,
  guarantees that no positive measure set of geodesics can be
  divergent by Poincar\'e recurrence.
\end{proof}

The construction of an infinite ray at $Y \in \Teich(S)$ asymptotic to
a given ray at $X \in \Teich(S)$ is an essential tool in our
discussion.  This is a general feature of complete $\CAT(0)$ spaces,
as shown in \cite[II.8, 8.3]{Bridson:Haefliger:npc}, and thus applies
to the completion $\closure{\Teich(S)}$.  More care is required,
however, to show that the resulting infinite ray in
$\closure{\Teich(S)}$ actually determines an infinite ray in
$\Teich(S)$.  Indeed, the possibility that a limit of unbounded or
even infinite geodesics might be finite cannot be ruled out a priori,
as was shown in \cite{Brock:nc} (see also \cite{Wolpert:compl}).
This is also a consequence of Proposition~\ref{prop:defined}.
Theorem~\ref{theorem:vis} follows from a key application of
Theorem~\ref{theorem:non-refraction}, the non-refraction of geodesics
in the Weil-Petersson completion.
 \begin{named}{Theorem~\ref{theorem:vis}}{\sc (Boundary at Infinity)} 
\vis 
\end{named}

\smallskip

\bold{Remark.}  Because of totally geodesic flats in the completion
arising from product strata, the condition that rays be merely
asymptotic, namely, that they remain a bounded distance apart, cannot
be improved to the condition that they be strongly asymptotic, though
we will see this follows for recurrent rays
(Theorem~\ref{theorem:recurrent:asymp}).

\begin{proof} It is a general consequence of \cite[II.8,
  8.3]{Bridson:Haefliger:npc} applied to the complete $\CAT(0)$ space
  $\closure{\Teich(S)}$ that we have a unique infinite geodesic ray
  $\ray'(t)$ in $\closure{\Teich(S)}$ based at $Y$ in the asymptote
  class of $\ray$ based at $X$.  Indeed, the ray $\ray'(t)$ is the
  limit of finite-length geodesics $\overline{Y\ray(t)}$ joining $Y$ to
  points along the ray $\ray$ with their parametrizations by
  arclength, a fact we note for future reference.  

  It remains only to conclude that $\ray'(t) \in \Teich(S)$ for each
  $t >0$.  But by Theorem~\ref{theorem:non-refraction}, for each $T>0$
  the geodesic $\ray'([0,T])$ has interior $\ray'((0,T))$ in the
  stratum $\calS_{\sigma_0\cap \sigma_T}$ where $\ray'(0) \in
  \calS_{\sigma_0}$ and $\ray'(T) \in \calS_{\sigma_T}$.  But since $Y
  \in \Teich(S)$ we have $\sigma_0 = \nullset$, so $\ray'(t)$ lies in
  the main stratum $ \calS_{\nullset} = \Teich(S)$ for each $t < T$.
  Since $T$ is arbitrary, the conclusion follows.

  It is general for a $\CAT(0)$ space that given a basepoint $X$, and
  an infinite ray $\ray$ at $X$, the ray $\ray$ is the unique
  representative of its asymptote class that is based at $X$.  Thus,
  we have a unique infinite ray based at a fixed $X$ in each asymptote
  class.  Applying the $\CAT(0)$-geometry of $\closure{\Teich(S)}$, it
  follows that if $\ray_n$ is a sequence of rays based at $X$ with
  convergent initial tangents to the initial tangent of the infinite
  ray $\ray_\infty$, then the corresponding infinite rays $\ray'_n$
  based at $Y$ in the same asymptote class converge to the ray
  $\ray'_\infty$ based at $Y$ in the same asymptote class as
  $\ray_\infty$.  Thus the change of basepoint map is a homeomorphism
  on the infinite rays.
\end{proof}

We remark that the assumption that $Y$ lies in the interior of
$\Teich(S)$ is just for simplicity: the same argument may be carried
out to prove the following stronger statement.

\begin{theorem}
  Let $\sigma$ and $\sigma'$ be simplices in $ \widehat{\calC(S)}$.
 Let $Y$ lie in the interior of a boundary stratum $\calS_\sigma.$
  Then given an infinite ray $\ray$ in $\closure{\Teich(S)}$ based at
  $X \in \calS_{\sigma'}$, there is a unique infinite ray $\ray'$ based at $Y$
  with $\ray'(t) \in \Teich(S)\cup \calS_\sigma$ for each $t$ so that
  $\ray'$ lies in the same asymptote class as $\ray$.
\label{theorem:visbar}
\end{theorem}

\begin{proof}
The proof goes through as before with the additional observation that for
each $s$ the limit $g_\infty([0,s))$ lies in $\Teich(S) \cup
\calS_\sigma$ by Theorem~\ref{theorem:non-refraction}.
\end{proof}

\section{Ending laminations and recurrent geodesics}
\label{section:elc}
The primary goal of this section is to establish
Theorem~\ref{theorem:relc}.
\begin{named}{Theorem~\ref{theorem:relc}}
{\sc (Recurrent Ending Lamination Theorem)}
\relc
\end{named}

The main technical tool in this section will be the following
application of the Gauss-Bonnet theorem.
\begin{theorem}
\label{theorem:recurrent:asymp}
  Let $\ray$ be a recurrent Weil-Petersson geodesic ray.  Then if $\ray'$ is
  a ray asymptotic to $\ray$ then $\ray$ is strongly asymptotic to $\ray'$.
\end{theorem}

We wish to harness the fact that the recurrent ray $\ray$ returns to a
portion of $\calM(S)$ where the sectional curvatures are definitely
bounded away from $0$.  To do this we employ the technique of
simplicial ruled surfaces, similar to Thurston's {\em
  pleated annulus} argument (cf. \cite{Thurston:hype2} and similar methods
in \cite{Canary:thesis} -- see also
\cite{Bonahon:tame,Canary:ends,Souto:cores}).

For the purposes of the proof we make the following definition:
\begin{defn}
  Given a Weil-Petersson geodesic ray $\ray \colon [0,T] \to
  \Teich(S)$ parametrized by arclength, and an $\epsilon >0$, we say a
  $t>0$ is an {\em $\epsilon$-recurrence} for $\ray$ if $\ray(t)$ is
  $\epsilon$-thick.  Given $\delta >0$, a collection $\{t_k\}$ of
  $\epsilon$-recurrences for a ray $\ray$ is {\em $\delta$-separated} if $|t_k -
  t_{k-1}| >\delta$.
\end{defn}

\begin{proof}[Proof of Theorem~\ref{theorem:recurrent:asymp}] Let
  $\ray$ and $\ray'$ be asymptotic rays, and fix a parametrization of
  $\ray$ for which $\ray(t)$ is the nearest point projection on $\ray$
  from the point $\ray'(t)$, where $\ray'(t)$ is parametrized by arclength.

Given points $X$, $Y$, and $Z$ in $\Teich(S)$, let $\triangle (X Y Z)$
denote the ruled triangle with vertices $X$, $Y$, and $Z$ ruled by
geodesics joining $Z$ to points along $\overline{X \; Y}$.
Then given $T>0$, we consider the ruled triangle 
$$\Delta_T  = 
\triangle (\ray(0) \ray(T) \ray'(0)).$$

The Weil-Petersson Riemannian metric induces a smooth metric
$\sigma_T$ on $\Delta_T$  whose Gauss curvature
is pointwise bounded from above by the upper bound on the
ambient sectional curvatures. 
\makefig{The ruled triangle
  $\Delta_T$.}{figure:simplicial}{\psfig{file=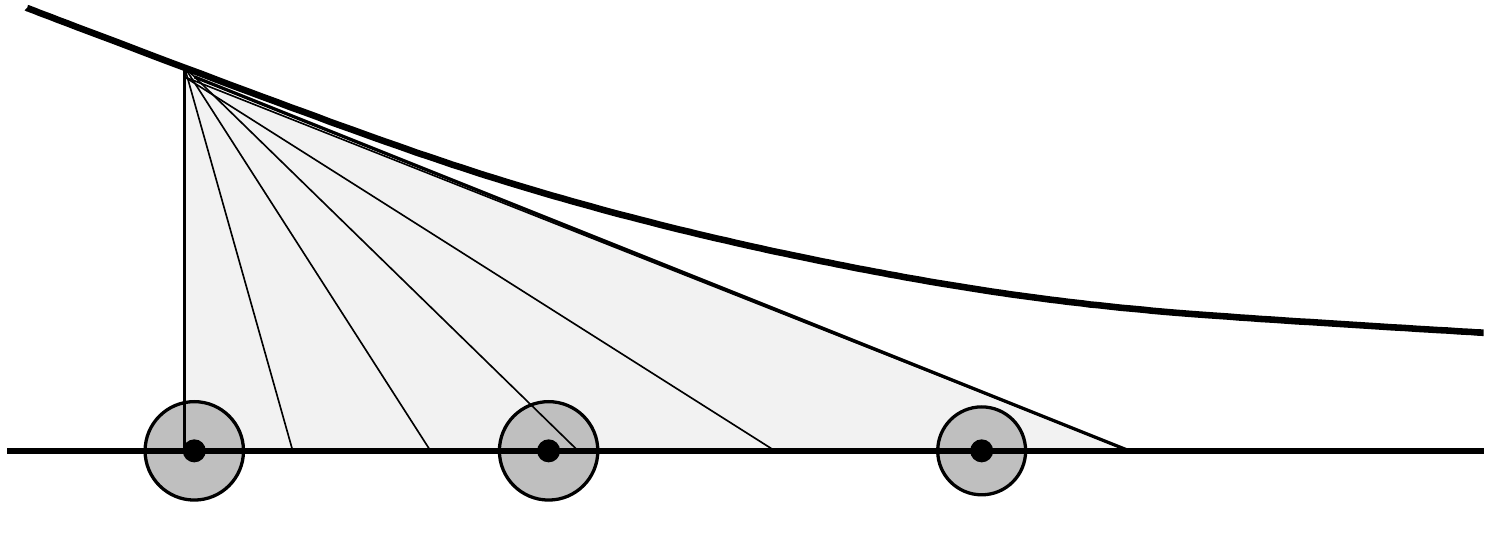,height=1.5in}}

By recurrence of $\ray$, there is an $\epsilon > 0$ so that for any
$\delta>0$ there is an infinite collection 
$\delta$-separated $\epsilon$-recurrences for $\ray$.  
Since the $\epsilon$-thick part projects to a compact subset of $\calM(S)$
\cite{Mumford:compact}, it follows that there is a $\delta_0>0$ so
that if $t$ is an $\epsilon$-recurrence for $\ray$ we have
\begin{equation}
\label{deltanaught}
\calN_{\delta_0}(\ray(t)) \subset \Teich_{\ge \epsilon/2}(S)
\end{equation}
  Thus, the sectional curvatures on
$\calN_{\delta_0}(\ray(t))$ are bounded from above by a negative
constant $\kappa_{\epsilon/2} < 0$, and thus the Gaussian curvature of
$\sigma_T$ on 
$$\Delta_T \cap
\calN_{\delta_0}(\ray(t))$$ is bounded above by $\kappa_{\epsilon/2}$
as well.

Assume there is a $\delta \in (0, \delta_0)$ so that the
distance from $\ray(t)$ to $\ray'(t)$ is bounded below by $\delta$.
The segments $\geod_T = \overline{\ray'(0)\; \ray(T)}$ converge to the
infinite ray $\ray'$, as $T \to \infty$, so for fixed $t$ we have
$\geod_T(t) \to \ray'(t)$ as $T \to \infty$.  It follows that for each
$t>0$ for which $\ray(t) \in \Teich_{\ge \epsilon}(S)$, there is a $T
>t$ so that $\geod_{T}(t)$ has distance at least $\delta/2$ from
$\ray(t)$.

We note that the inclusion map on $\Delta_T$ is $1$-Lipschitz from the
$\sigma_T$-metric to the Weil-Petersson metric.  
Then the
intersection of the triangle $\Delta_T$ with 
the neighborhood $\calN_{\delta/2}(\ray(t))$ contains a region with area
at least $\pi {\delta}^2/16$ in the intrinsic metric on $\Delta_{T}$, since
each such intersection contains a sector in $\Delta_{T}$ with radius $\delta/2$
and angle $\pi/2$ whose $\sigma_T$-area is minorized by the area of a
Euclidean sector with the same radius and angle.  

Let $\{ t_k \}_k$ be an infinite collection of $\delta$-separated
$\epsilon$-recurrences for $\ray$. Given $N >0$ take 
$T_N$ so that $\geod_{T_N}(t_k)$ has distance at least $\delta/2$ from
$\ray(t_k)$ for each $k \le N$.
Then we have the estimate
\begin{equation}
\label{definite}
\left|\int_{\Delta_{T_N}}\kappa \; dA\right| > N \left|\kappa_{\epsilon/2}\right| \frac{\pi
  {\delta}^2}{16}
\end{equation}
on the absolute value of the integral of the Gauss curvature $\kappa$ over
$\Delta_{T_N}$.  

The Gauss-Bonnet Theorem bounds the integral (\ref{definite}) from
above by $\pi$, independent of $T_N$, so letting
\begin{equation}
\label{recurrences}
N(\epsilon,\delta) =  \frac{16}{ \left|\kappa_{\epsilon/2}\right| {\delta}^2}
\end{equation}
the bound $N < N(\epsilon,\delta)$
to the cardinality  of the set $\{t_k\}_{k=0}^N$ follows.  

Since $N$ is arbitrary, (\ref{recurrences}) contradicts the recurrence
of $\ray$ to the $\epsilon$-thick part, and we conclude the existence
of $\hat{t}$ for which $\geod_T(\hat{t}) \cap
\calN_{\delta/2}(\ray(\hat{t})) \not= \emptyset$ independent of $T$.
Since we have $\geod_T(\hat{t}) \to \ray'(\hat{t})$ as $T \to \infty$,
we conclude that $\ray'(\hat{t}) \cap \calN_{\delta}(\ray(\hat{t}))
\not= \emptyset$.  Since $\delta>0$ is arbitrary, and since the
distance from $\ray(t)$ to $\ray'(t)$ is a non-increasing function in
a $\CAT(0)$ space, it follows that the rays are strongly asymptotic.
\end{proof}

\bold{Remark:} M. Bestvina
and K. Fujiwara have observed indepenently the applicability of this
ruled surface technique to the study of action of $\Mod(S)$
on $\closure{\Teich(S)}$ as the isometry group of a $\CAT(0)$-space
(cf. \cite{Bestvina:Fujiwara:rank}).
\smallskip

We employ the fact that recurrent rays exhibit such strongly
asymptotic behavior to conclude Theorem~\ref{theorem:rvis}.

\begin{named}{Theorem~\ref{theorem:rvis}}
{\sc (Recurrent Visibility)}
\rvis
\end{named}

\begin{proof}
  We seek to exhibit a bi-infinite Weil-Petersson geodesic $\geod \colon
  \reals \to \Teich(S)$ with the property that $\geod$ is strongly
  asymptotic to the recurrent ray $\ray^+$ in positive time and
  asymptotic to $\ray^-$ in negative time.  In other words, we claim there is a
  reparametrization $t \mapsto s(t) >0$ so that we have
$$d(\geod(s(t)),\ray^+(t))  \to 0$$ as $t \to \infty$, and 
  $$d(\geod(t),\ray^-(-t))$$ is bounded for $t < 0$.

  Consider geodesic chords $\geod_n$ joining $\ray^+(n)$ to
  $\ray^-(n)$ for integers $n >0$.  We wish to show that after passing
  to a subsequence there are paramater values $\hat{t}_n$ so that
  $\geod_n(\hat{t}_n)$ converges to some $ Z \in \Teich(S)$.

  Assume $\ray^+$ is recurrent to the $\epsilon$-thick part.  For this
  $\epsilon$, let $\delta_0$ be chosen as in (\ref{deltanaught}).
  Then given $\delta \in (0, \delta_0)$ we let
  $\{t_k\}_{k=0}^\infty$ be $\delta$-separated $\epsilon$-recurrences
  for $\ray^+$.

  Since the rays $\ray^-$ and $\ray^+$ are assumed distinct, we may
  omit finitely many recurrences and assume that for each $k$ the
  distance of $\ray^+(t_k)$ from $\ray^-$ is at least $\delta$.
  Otherwise, the two rays lie in the same asymptote class, and are
  thus identical since they are based at the same point in a $\CAT(0)$
  space.

We claim that there is a $\hat{k}$ so that 
$$\geod_n \cap \calN_{\delta/2}(\ray^+(t_{\hat{k}})) \not= \emptyset$$ for all $n$.
Let $s_0 >0$ be the minimal parameter so that $\ray^-(s_0)$ has distance
at least $\delta$ from $\ray^+$.  Then for $n>s_0$, consider the ruled triangle
in $\Teich(S)$ given by
$$\calT_n = \triangle (\ray^+(0)\ray^+(n)\ray^-(n)).$$ 
If $N$ is chosen so that  $t_N < n$, and
$$\calN_{\delta/2}(\ray^+(t_N)) \cap \geod_n = \emptyset,$$
then as in (\ref{recurrences}) we have the upper
bound
$$N < N(\epsilon, \delta).$$

Since $\{t_k\}_k$ is infinite, we conclude there is a $\hat{k}$ so
that $\geod_n$ intersects $\calN_{\delta/2}(t_{\hat{k}})$ for all $n$
sufficiently large.  As the $\delta/2$-ball
$\calN_{\delta/2}(\ray^+(t_{\hat{k}}))$ lies in $\Teich_{\ge
  \epsilon/2}(S)$, it is precompact, and we may pass to a convergent
subsequence so that points $Z_n \in \geod_n$ converge to a limit $Z$
in $\Teich(S)$.

Reparametrizing $\geod_n$ by arclength so that $\geod_n(0) = Z_n$, we
observe that for each $t \in \reals^+$ the sequence $\{\geod_n(t)\}_n$ is
Cauchy.  To see this, note that if $\mathbf{h}_n^+$ is the segment joining
$Z$ to $\ray^+(n)$ parametrized by arclength, choosing $n_t$ so that 
so that $d(Z, \ray^+(n) ) >t$  for each $n \ge n_t$, the sequence
$\{\mathbf{h}_n^+(t)\}_{n = n_t}^\infty$ is Cauchy  
(as in \cite[II.8,8.3]{Bridson:Haefliger:npc} and its use in
Theorem~\ref{theorem:vis}). The assertion then follows from the
observation that $\CAT(0)$ geometry
guarantees
$$d(\mathbf{h}_n^+(t),\geod_n^+(t)) < d(Z_n, Z) \to 0.$$ 
The symmetric argument shows $\{\geod_n(-t)\}_n$ is Cauchy as well.

There is thus a limiting geodesic $\geod \colon \reals \to
\overline{\Teich(S)}$ for which the ray $\geod^+ = \geod
\vert_{[0,\infty)}$ lies in the asymptote class of $\ray^+$ and
$\geod^- = \geod \vert_{(-\infty,0]}$ lies in the asymptote class of
$\ray^-$.  Since $\geod^+$ and $\geod^-$ are the unique
representatives of the asymptote classes of $\ray^+$ and $\ray^-$
based at $Z$, an application of Theorem~\ref{theorem:vis} ensures
$\geod^+$ and $\geod^-$ lie entirely within $\Teich(S)$ as claimed.

For statement~(2), we note that by Theorem~\ref{theorem:vis}, the ray
$\geod^+$ is the limit of geodesics $\geod_n^+ $ joining $\geod(0) =
Z$ to $\ray^+(n)$, so if $\mu$ has bounded length along $\ray^+$, then
convexity of the length of $\mu$ guarantees that the length of $\mu$
is uniformly bounded on $\geod^+$.  The same argument applies to
$\geod^-$.  Statement~(2) follows.
\end{proof}

In Section~\ref{section:preliminaries}, we employed the boundedness of
length functions for ending measures along a ray to establish that the
ending lamination is well defined.  For a recurrent ray, however, we
can guarantee that the length of any lamination with bounded length
decays to zero.
\begin{lem}
Let $\ray(t)$ be a recurrent ray, and let $\mu \in \ml(S)$ be any
lamination with $\ell_\mu({\ray(t)}) <K$ along $\ray(t)$.   Then we have
$$\ell_\mu({\ray(t)}) \to 0$$ as $t \to\infty$.
\label{lemma:recurrent:pinches}
\end{lem}

\begin{proof}
  Assume $\ray(t)$ recurs to the $\epsilon$-thick part at times $t_n
  \to \infty$.  Wolpert's extension of his convexity theorem for
  geodesic length functions guarantees that the length of $\mu \in
  \ml(S)$, in addition to being convex along geodesics
  \cite{Wolpert:Nielsen}, satisfies the following stronger convexity
  property: given $\epsilon >0$, there is a $c >0$ so that at each $t$
  for which $\ray(t)$ lies in the $\epsilon$-thick part, we have
  $$\ell_\mu''({\ray(t)}) > c \ell_\mu({\ray(t)})$$ (see
  \cite{Wolpert:glf}).  The proof of the Lemma then follows from the
  observation that if the bounded convex function
  $\ell_\mu({\ray(t)})$ does not tend to zero, then we nevertheless
  have $\ell_\mu({\ray(t)}) \to C >0$ as $t\to \infty$, which guarantees
  that $\ell_\mu''({\ray(t)}) \to 0$ by convexity.  This contradicts
  the above inequality at the times $t_n$ for $n$ sufficiently large.
\end{proof}

Though the ending lamination need not fill the surface in general, the
recurrent rays provide a class of rays where each lamination with
bounded length along the ray fills $S$.
\begin{prop}
  \label{prop:mu:filling}
Let $\mu$ be any measured lamination with bounded length along the
recurrent ray $\ray(t)$.  Then $\mu$ is a filling lamination.
\end{prop}

\begin{proof}
  Assume $\mu$ does not fill, and let $S(\mu)$ be the supporting
  subsurface for its support $|\mu|$.  Let $\gamma_n \in
  \calC(S(\mu))$ be a sequence of simple closed curves whose
  projective classes $[\gamma_n]$ converge to $[\mu]$ in $\pl(S)$.
  Note in particular that $$i(\bdry S(\mu),\gamma_n) = 0$$ for each
  $n$.

We claim that 
given any $Z \in \Teich(S(\mu))$ there is a Weil-Petersson ray $\hat
  \ray$ in $\Teich(S(\mu))$ based at $Z$ along which $\mu$ has bounded
  length.
  To see this, note that any limit $\hat \ray$ of finite-length rays
  $\overline{Z \; Z_n}$ joining $Z$ to nodal surfaces $Z_n$ in
  $\closure{\calS_{\gamma_n}} \cap \closure{\Teich(S(\mu))}$ has the
  property that $[\mu]$ is the projective class of a lamination with
  bounded length along $\hat \ray$ by Lemma~\ref{lem:continuity}.  The
  fact that $\mu$ fills $S(\mu)$ guarantees that $\hat \ray$ has no
  pinching curves.  Thus $\hat \ray$ has infinite length.

Letting $\sigma_\mu \in \widehat{\calC(S)}$ be the simplex spanned by
the curves in $\bdry S(\mu)$  we note that the stratum
$\calS_{\sigma_\mu}$ is the metric product of Weil-Petersson metrics
on $\Teich(S(\mu))$ and the Weil-Petersson metrics on $\Teich(Y)$
where $Y$ is the disjoint union of  non-annular components of
$S\setminus S(\mu)$.

Together with the basepoint $X$, then, the ray $\hat \ray$ naturally
determines a ray $\closure{\ray}$ in the stratum $\calS_{\sigma_\mu}$
by taking the projection of $\closure{\ray(t)}$ to $\Teich(S(\mu))$ to
be $\hat \ray (t)$ and identifying each other coordinate of
$\closure{\ray}(t)$ in the product decomposition of
$\calS_{\sigma_\mu}$ with the (constant) coordinate function of the
nearest point projection of $X$ to $\calS_{\sigma_\mu}$.

Applying Theorem~\ref{theorem:visbar}, there is a unique ray $\ray'$
based at $X$ asymptotic to $\closure{\ray}$.  The ray $\ray'$ is
constructed as a limit of segments $\geod_t = \overline{X \; \closure{\ray}(t)}$ joining
$X$ to points along $\closure{\ray}$.  The length of $\mu$ and each
curve $\gamma \subset \bdry S(\mu)$ is uniformly bounded on the segments
$\overline{X \; \closure{\ray}(t)}$, by convexity of length functions.
Applying continuity of length, then, we have a $K>0$ so that
$$\ell_\mu({\ray'(t)}) < K \
 \ \ \text{and} \ \ \ \ell_\gamma({\ray'(t)}) < K$$ for each
$\gamma \in \bdry S(\mu)$.

If $\ray'$ is distinct from $\ray$, however,
Theorem~\ref{theorem:rvis} guarantees that we may find a bi-infinite
geodesic $\geod$ whose forward trajectory is strongly asymptotic to
$\ray$, by recurrence, and so that $\geod \vert_{(-\infty,0]}$ stays a
bounded distance from $\ray'$.  Once again, the length
$\ell_\mu({\geod(t)})$ of $\mu$ is uniformly bounded over the entire
bi-infinite geodesic $\geod$.  Since $\ell_\mu({\geod(t)})$ approaches
$0$ as $t \to \infty$ by Lemma~\ref{lemma:recurrent:pinches}, the
boundedness of $\ell_\mu({\geod(t)})$ yields a contradiction to its
strict convexity.  We conclude that $\ray = \ray'$ and thus that
$\gamma$ has bounded length along the ray $\ray(t)$.  But applying
Lemma~\ref{lemma:recurrent:pinches} once again, boundedness implies
that the length of $\gamma$ tends to zero along $\ray(t)$, violating
recurrence of $\ray(t)$.

We conclude that $\mu$ fills $S$.
\end{proof}

Theorem~\ref{theorem:relc} will be a direct consequence of the
following characterization of measures that arise with bounded length
along a recurrent Weil-Petersson ray, together with Theorem~\ref{theorem:rvis}.
\begin{lem}
\label{lemma:bounded}
  Let $\ray$ be a recurrent ray with ending lamination $\elam =
  \elam(\ray)$.  If $\mu \in \ml(S)$ then the following are
  equivalent:
\begin{enumerate}
\item \label{item:support} The lamination $\mu$ has support $\elam$.
\item \label{item:bounded} The length function $\ell_\mu(\ray(t))$ is bounded.
\item \label{item:decays} $\ell_\mu({\ray(t)}) \to 0.$
\end{enumerate}
\end{lem}

\begin{proof}
  We first verify that (\ref{item:support}) implies
  (\ref{item:bounded}).  Let $\Sigma$ denote the simplex of projective
  classes of measures on $\elam$ in $\pl(S)$, and let $\hat{\mu} \in
  \ml(S)$ be a representative of the projective class determined by a
  point in the interior of the top dimensional face.  Then $\hat{\mu}$
  is a positive linear combination of all ergodic measures on $\elam$.

  Let $\gamma_n$ be a sequence of simple closed curves for which the
  projective classes $[\gamma_n]$ converge to $[\hat{\mu}]$, and let
  $\ray_n$ be a sequence of finite rays based at $X$ limiting to
  points $Z_n$ in the strata $\calS_{\gamma_n}$.  Since $\gamma_n$ are
  pinching curves for $\ray_n$, Lemma~\ref{lem:continuity} guarantees
  that any limit $\ray_\infty$ of a convergent subsequence of $\ray_n$
  has the property that $\hat{\mu}$ has bounded length along $\ray_\infty$.
  Since $\hat{\mu}$ is a positive linear combination of all the ergodic
  measures on $\elam$, it follows that each ergodic measure on $\elam$
  has bounded length along $\ray_\infty$.  Hence, any measured
  lamination representing a projective class in $\Sigma$ has bounded
  length along $\ray_\infty$ since each is a linear combination of
  ergodic measures.

 Since $\elam(\ray)$ is
  filling, by Proposition~\ref{prop:mu:filling}, we have that $\hat{\mu}$ is
  filling.  This guarantees that $\ray_\infty$ has infinite length,
  since otherwise $\ray_\infty$ would have a pinching curve $\gamma$
  with $i(\gamma,\hat{\mu}) >0$, violating the length bound on $\hat{\mu}$ along
  $\ray_\infty$.

  If $\bar{\mu}$ is any ending measure for $\ray$, then $\bar{\mu}$
  represents a projective class in $\Sigma$, and thus has bounded
  length along $\ray_\infty$.  If $\ray$ and $\ray_\infty$ are
  distinct rays, then Theorem~\ref{theorem:rvis} guarantees that we
  have a bi-infinite geodesic $\geod(t)$ asymptotic to $\ray$ and
  $\ray_\infty$ along which $\bar{\mu}$ has bounded length, which
  contradicts strict convexity of the length function for $\bar{\mu}$ along
  $\geod$.  It follows that $\ray = \ray_\infty$. 

  Since $\mu$ also represents a measure in $\Sigma$, it follows that
  $\mu$ has bounded length along $\ray$, verifying that
  (\ref{item:support}) implies (\ref{item:bounded}).

  Conclusion (\ref{item:decays}) follows from conclusion
  (\ref{item:bounded}) by an application of
  Lemma~\ref{lemma:recurrent:pinches}.

  By Corollary~\ref{cor:bounded:length}, (\ref{item:decays}) implies
  that any ending measure $\mu'$ for $\ray$ has the property that
$$i(\mu,\mu') = 0.$$
By Proposition~\ref{prop:mu:filling} each of $\mu$ and $\mu'$ fills,
so they have identical support, verifying conclusion
(\ref{item:support}), and hence proving the Lemma.
\end{proof}

We are now ready to prove Theorem~\ref{theorem:relc}, that recurrent
rays with the same ending lamination are strongly asymptotic.

\begin{proof}[Proof of Theorem~\ref{theorem:relc}]
  Let $\ray$ be a recurrent ray based at $X \in \Teich(S)$, with
  ending lamination $\elam = \elam(\ray)$.  
  Let $\ray'$ be another ray based at $X \in \Teich(S)$ with ending
  lamination $\elam$, and let $\mu$ be an ending measure for $\ray'$.
  Then $\mu$ has support $\elam$, so by an application of
  Lemma~\ref{lemma:bounded} $\mu$ has bounded length along $\ray$ as
  well.

  If $\ray$ and $\ray'$ are distinct rays, then
  Theorem~\ref{theorem:rvis} guarantees that we have a bi-infinite
  geodesic $\geod(t)$ asymptotic to $\ray$ and $\ray'$ along which
  $\mu$ has bounded length, which contradicts strict convexity of the
  length function for $\mu$ along $\geod$.  It follows that $\ray =
  \ray'$.

  If $\ray''$ is a ray based at $Y \not= X$, with ending lamination
  $\elam$, there is a unique ray $\ray'$ based at $X$ in the asymptote
  class of $\ray''$ by Theorem~\ref{theorem:vis}.  Applying the above
  argument to $\ray'$ we may conclude that $\ray''$ and $\ray$ are in
  the same asymptote class.
  Theorem~\ref{theorem:recurrent:asymp} then guarantees that
  $\ray''$ and $\ray$ are strongly asymptotic, concluding the proof.
\end{proof}

As a further consequence, we note the following.  
\begin{cor}
  Let $\geod(t)$ be a bi-infinite Weil-Petersson geodesic whose
  forward trajectory is recurrent.  Then the ending laminations
  $\elam^+$ and $\elam^-$ for the rays $\geod^+ = \{\geod(t)\}_{t = 0}^\infty$
  and $\geod^- = \{\geod(t)\}_{t= 0}^{-\infty}$ bind the surface $S$.
\label{cor:binding}
\end{cor}

\begin{proof}
  The ending lamination $\elam^+$ for the forward trajectory fills the
  surface, so the ending lamination for the backward trajectory must
  intersect it, since otherwise the laminations $\elam^-$ and
  $\elam^+$ would be identical therefore we would have $\geod^+ =
  \geod^-$ by Theorem~\ref{theorem:relc}, a contradiction.
\end{proof}

To derive Corollary~\ref{cor:parametrize}, we establish a final further
continuity property for ending measures when the limit is recurrent.
\begin{prop}
  If $\ray_n$ is a convergent sequence of rays at $X$ with a recurrent
  limit $\ray$, any sequence $\mu_n$ of ending measures or pinching
  curves for $\ray_n$ converges in $\pl(S)$ up to subsequence to a
  measure on $\elam(\ray)$.
\label{prop:continuity:recurrent}
\end{prop}

\begin{proof}
Let $\mu$ be any limit of $\mu_n$ in $\pl(S)$ after passing to a
subsequence.  Then by Lemma~\ref{lem:continuity}, the length
$\ell_\mu({\ray(t)})$ is bounded.  Since $\ray$ is recurrent, any
ending measure $\mu'$ for $\ray$ 
fills $S$ by Proposition~\ref{prop:mu:filling}.
But by
Corollary~\ref{cor:bounded:length}, we have 
$$i(\mu,\mu') = 0$$  so $\mu$ and $\mu'$ have identical support since
$\mu'$ is filling.  Hence, $\mu$ is a measure on $\elam(\ray)$.
\end{proof}

Restricting to the recurrent rays, we obtain Corollary~\ref{cor:parametrize}.
\begin{named}{Corollary~\ref{cor:parametrize}}
\parametrize
\end{named}

\begin{proof}
  That the map is a bijection follows from the fact that $\rel(S)$ is
  defined as its image and from Theorem~\ref{theorem:relc}.  

  To show continuity in each direction, we begin by noting that
  although the topology induced by forgetting the measure on a measured
  lamination is not a Hausdorff topology on the geodesic laminations
  admitting measures, it is Hausdorff when one restricts to those that
  fill the surface, namely, the subset $\el(S)$ (see \cite[\S 7]{Klarreich:boundary}).  As such it suffices
  to consider sequential limits to establish continuity.  

  Furthermore, the topology of convergence of asymptote classes of of
  infinite rays in $\bdry_\infty \closure{\Teich(S)}$ agrees with the
  topology of convergence in the visual sphere of representatives
  emanating from a given basepoint $X$ as in Theorem~\ref{theorem:vis}
  (see also \cite[II.8,8.8]{Bridson:Haefliger:npc}).

  Let $\ray_n$ be a sequence of recurrent rays based at $X$ with recurrent limit
  $\ray$ based at $X$.  By Proposition~\ref{prop:mu:filling}, their ending
  laminations $\elam_n$ are filling laminations and thus determine
  points in $\rel(S)$.  Their recurrent limit $\ray$ has ending
  lamination $\elam(\ray)$, with support identified with the support
  of a limit of measures on $\elam_n$ by
  Proposition~\ref{prop:continuity:recurrent}, so $\elam$ is the limit
  of $\elam_n$ in $\rel(S)$, by the definition of the topology on
  $\el(S)$.

  For continuity in the other direction, compactness of the visual
  sphere at $X$ guarantees that any sequence of laminations
  $\lambda_n$ converging to $\lambda_\infty$ in $\rel(S)$ determine a
  sequence of rays $\ray_n$ at $X$ with limit $\ray_\infty$ at $X$
  after passing to a subsequence.  A convergent family of measures
  $\mu_n$ on $\lambda_n$ has limit $\mu_\infty$, a measure on
  $\lambda_\infty$, with bounded length on the limiting ray
  $\ray_\infty$ by Lemma~\ref{lem:continuity}.  Since $\mu_\infty$ is
  filling, and any ending measure or weighted pinching curve $\mu'$
  for $\ray_\infty$ satisfies $i(\mu,\mu') = 0$, we conclude that
  $\mu'$ has the same support as $\mu_\infty$, namely
  $\lambda_\infty$.  Thus $\ray_\infty$ is the recurrent ray at $X$
  determined (uniquely) by $\lambda_\infty$.  Since any accumulation
  point of the rays $\ray_n$ has this property, the original sequence
  of rays itself was convergent to $\ray_\infty$, obviating passage to
  subsequences.
\end{proof}

\section{The topological dynamics of the geodesic flow}
\label{section:dynamics}
We now relate the preceding results to the study of the Weil-Petersson
geodesic flow on $\calM^1(S)$.

Though it is seen in \cite{Brock:nc} that the change of basepoint map
is discontinuous on the visual sphere, the visibility property for
recurrent rays (Theorem~\ref{theorem:rvis}) is sufficient to remedy
the situation for considerations of topological dynamics, yielding
Theorem~\ref{theorem:periodic}, whose proof we now supply.

\begin{named}{Theorem~\ref{theorem:periodic}}{\sc (Closed Orbits Dense)}
  The set of closed Weil-Petersson geodesics is dense in $\calM^1(S)$.
\end{named}

\begin{proof}
  Because of the density of doubly recurrent geodesics in the unit
  tangent bundle $T^1 \Teich(S)$, it suffices by a diagonal argument
  to approximate a bi-recurrent direction with periodic geodesics.

  Let $\{\geod(t)\}_{t = -\infty}^\infty$ be a bi-infinite geodesic
  that is doubly recurrent.  Let $X = \geod(0)$ be a basepoint, and
  let $\elam^+$ be the ending lamination for the forward ray
  $\geod^+(t) = \{\geod(t) \}_{t = 0}^{\infty}$ and likewise let
  $\elam^-$ denote the ending lamination for the backward ray
  $\geod^-= \{\geod(-t)\}_{t = 0 }^{\infty}$.  By
  Corollary~\ref{cor:binding}, $\elam^+$ and $\elam^-$ bind the
  surface $S$, so letting $\mu^+$ and $\mu^-$ be measures on $\elam^+$
  and $\elam^-$, respectively, any pair of simple closed curves
  $\gamma^+$ and $\gamma^-$ very close to $\mu^+$ and $\mu^-$ in
  $\pl(S)$ also bind $S$.

Letting $\tau_+$  be a Dehn twist about $\gamma^+$ and $\tau_-$ be a
Dehn twist about $\gamma^-$, the composition
$$\psi_k = \tau_+^k \compos \tau_-^{k}$$ is pseudo-Anosov for all $k$
sufficiently large \cite{Thurston:mcg}.  Furthermore, the stable and
unstable laminations for $\psi_k$ converge to $\gamma^+$ and
$\gamma^-$ in $\pl(S)$ as $k \to \infty$.  Diagonalizing, then, we
obtain a sequence of pseudo-Anosov mapping classes $\varphi_n$ whose
unstable and stable laminations $\mu_n^+$ and $\mu_n^-$ converge to
$\mu^+$ and $\mu^-$ in $\pl(S)$.  Since the supports $|\mu_n^\pm|$ and
$|\mu^\pm|$ lie in $\rel(S)$, we have convergence of $|\mu_n^\pm|$ to
$\elam^\pm$ in $\rel(S)$ by the definition of the topology on
$\el(S)$.

Letting $A_n$ be the axis for $\varphi_n$, we claim $A_n$ is
arbitrarily close to $\geod$ at $\geod(0)$ in the unit tangent bundle for
$n$ sufficiently large.

To see this, we apply Theorem~\ref{theorem:vis} to obtain a ray
$\ray_n^+$ in $\calV_X(S)$ asymptotic to $A_n$ in the forward
direction.  We note that, as $A_n$ is itself doubly recurrent, the ray
$\ray_n^+$ is strongly asymptotic to $A_n$, by
Theorem~\ref{theorem:recurrent:asymp}, and that the ending lamination
$\elam^+_n$ for $\ray^+_n$ is equal to the support of $\mu_n^+$.  It
follows that $\elam^+_n$ converges to $\elam^+$ in $\rel(S)$.
Likewise, if $\ray_n^-$ denotes the ray in $\calV_X(S)$ asymptotic to
$A_n$ in the negative direction, then $\elam^-_n = \elam(\ray_n^-)$
converges to $\elam^-$ in $\rel(S)$.
The parametrization of recurrent rays by their ending laminations in
$\el(S)$, Corollary~\ref{cor:parametrize}, guarantees that $\ray_n^+$
and $\ray_n^-$ converge to $\geod^+$ and $\geod^-$ respectively.

Let $\epsilon>0$ be taken so that $\geod^+$ and $\geod^-$ recur to
$\Teich_{\ge 4\epsilon}(S)$.  
Assume $\delta_0>0$ is chosen as in
(\ref{deltanaught}) so that given $Z \in \Teich_{\ge \epsilon}(S)$ and
$\delta< \delta_0$ the $\delta$-neighborhood $\calN_\delta(Z)$ is
precompact in $\Teich(S)$.  Let $\{t_k\}_k$ and $\{s_k\}_k$ be
$4\delta$-separated $4\epsilon$-recurrences for $\geod^+$ and
$\geod^-$ respectively.

The convergence of $\ray_n^+ \to \geod^+$ and $\ray_n^- \to \geod^-$
guarantees that given any positive integer $N$, there is a $n_N$ so
that for $n>n_N$ the parameters $\{t_k\}_{k=0}^N$ and
$\{s_k\}_{k=0}^N$ are $2\delta$-separated $2\epsilon$-recurrences for
$\ray_n^+$ and $\ray_n^-$ respectively.  Letting
$$\geod_{n,T}^+ = \overline{X\; Z_n^+(T)} \ \ \ \text{and} \ \ \ 
\geod_{n,T}^- = \overline{X\; Z_n^-(T)}$$ be the geodesic segments
joining $X$ to nearest points $Z_n^+(T)$ and $Z_n^-(T)$ on
$A_n$ to $\ray_n^+(T)$ and $\ray_n^-(T)$, we have the
convergence of $\geod_{n,T}^+$ to $\ray_n^+$ and $\geod_{n,T}^-$ to
$\ray_n^-$ as $T \to \infty$ for fixed $n$.  There is a $T_N$, then, so that
$t_k$ and $s_k$ are $\delta$-separated $\epsilon$-recurrences
for $\geod_{n,T_N}^+$ and $\geod_{n,T_N}^-$ for each $k \le N$ and $n> n_N$.

\makefig{The axis $A_n$ converges to
  $\geod$.}{figure:ruled2}{\psfig{file=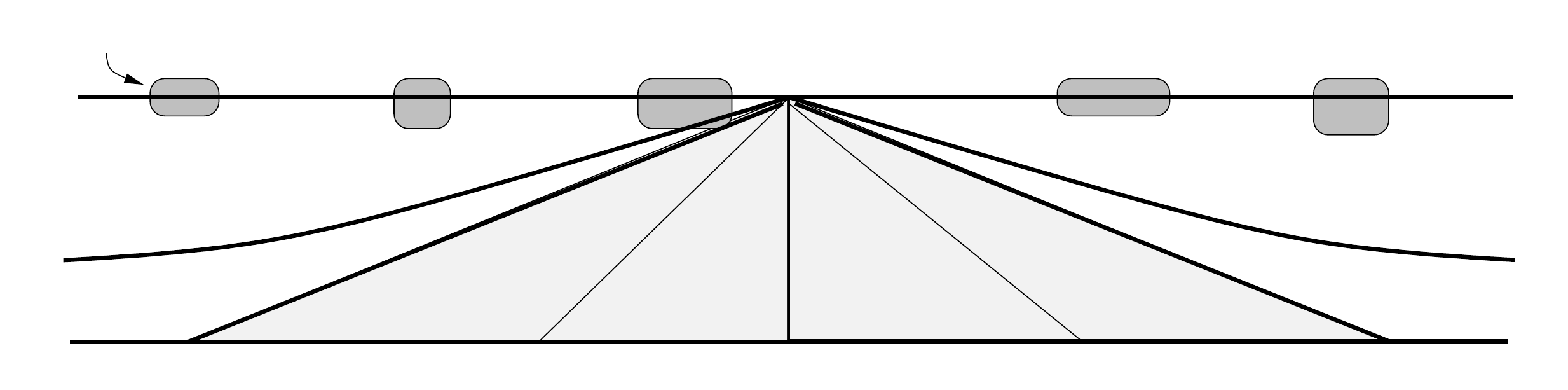,height=1.45in}}

We wish to apply the ruled triangle argument of 
Theorem~\ref{theorem:recurrent:asymp} to the 
two triangles $\Delta_{n}^+(T)$ and $\Delta_n^-(T)$ where
$$\Delta_n^+(T) = 
\triangle(Z_n^+(T)\; Z_n^0\; X) \ \ \ \text{and} \ \ \ 
\Delta_n^-(T) = 
\triangle(Z_n^-(T)\;
Z_n^0\; X)$$ 
(see figure~\ref{figure:ruled2}).

Taking $\delta < \delta_0$ and letting $N(\epsilon,\delta)$ be as in
the application of Gauss-Bonnet in (\ref{recurrences}) assume $n
> n_{N(\epsilon,\delta)}$ and $T > T_{N(\epsilon,\delta)}$.
If $N'$ is the maximal integer such that
$$\cup_{k=1}^{N'} 
\left( \calN_{\delta/2}(\geod_{n,T}^+(t_k))\right) \cap A_n = \emptyset \ \ \ \ \text{or} \
\ \ \ 
\cup_{k=1}^{N'}  
\left(\calN_{\delta/2}(\geod_{n,T}^-(s_k))\right) \cap A_n = \emptyset$$ 
it follows that 
\begin{equation*}
N' < N(\epsilon,\delta).
\end{equation*}
Taking $\hat{k} \ge N(\epsilon,\delta)$, then, we may conclude that for
$n$ and $T$ sufficiently large we have
$$\calN_{\delta/2}(\geod_{n,T}^+(t_{\hat{k}})) \cap A_n \not= \emptyset  \ \ \ 
\text{and}  \ \ \
\calN_{\delta/2}(\geod_{n,T}^-(s_{\hat{k}})) \cap A_n \not=
\emptyset.$$
Since for fixed $n$ we have $\geod_{n,T}^+(t_{\hat{k}}) \to
\ray_n^+(t_{\hat{k}})$ and $\geod_{n,T}^-(t_{\hat{k}}) \to
\ray_n^-(t_{\hat{k}})$ as $T \to \infty$, 
it follows that 
$$\calN_{\delta}(\ray_n^+(t_{\hat{k}})) \cap A_n \not= \emptyset  \ \ \ 
\text{and}  \ \ \
\calN_{\delta}(\ray_n^-(s_{\hat{k}})) \cap A_n \not= \emptyset$$
for $n$ sufficiently large.

Properties of $\CAT(0)$ geometry guarantee that the points
$Z^+_n(t_{\hat{k}})$ and $Z^-_n(s_{\hat{k}})$ on $A_n$ 
bound geodesic subsegments $\ell_{n,\hat{k}}$ of $A_n$ that converge
up to subsequence to a geodesic segment $\ell_{\hat{k}}$ within
$\delta$ of $\geod$ as $n \to \infty$.  Thus, a diagonal argument
allows us to extract a bi-infinite geodesic limit $A_\infty$ of the
axes $A_n$ so that for each $\delta>0$, and $T'>0$, $A_\infty$
contains a subsegment of width at least $T'$ within $\delta$ of
$\geod$.  We conclude $A_\infty = \geod$.

The projections of $A_n$ to $\calM(S)$ are closed
geodesics approximating the doubly recurrent projection of $\geod$ to
$\calM(S)$, as was desired.
\end{proof}
Using the boundary theory for the recurrent rays and its connection
with measured laminations, we can harness the north-south dynamics of
pseudo-Anosov elements on $\pl(S)$ to establish
Theorem~\ref{theorem:dense} as a consequence of
Theorem~\ref{theorem:periodic}.

\begin{named}{Theorem~\ref{theorem:dense}}{\sc (Dense Geodesics)}
  Given any $X \in \Teich(S)$, there is a Weil-Petersson geodesic ray
  based at $X$ whose projection to $\calM^1(S)$ is dense.
\end{named}

\begin{proof}
Given a pseudo-Anosov mapping class $\psi$ let $\mu^+$ and $\mu^-$
denote representatives of the attracting and repelling fixed points
$[\mu^+]$ and $[\mu^-]$ for the action of $\psi$ on $\pl(S)$,
respectively.  Let $A_\psi$ denote its axis in the Weil-Petersson
metric, and let $\lambda^+ = |\mu^+|$ and $\lambda^- = |\mu^-|$
denote the support in $\el(S)$ of the attracting and repelling
laminations for $\psi$.

Given $X \in \Teich(S)$, and $\delta>0$, we have from
Corollary~\ref{cor:parametrize} that there is a neighborhood
$U_\delta^+(\psi) \subset \rel(S)$ of $\lambda^+$ so that if $\lambda'
\in U_\delta^+(\psi)$ is the support of the attracting fixed point
$[\mu'] \in \pl(S)$ of another pseudo-Anosov element
$\psi'$, then there is a $Z \in A_\psi$ so that the ray $\ray$ based at $X$  with
ending lamination $\lambda'$ contains a segment within Hausdorff
distance $\delta$ of a full period
$g = \overline{Z \; \psi(Z)}$ of the action of $\psi$ on $A_\psi$.

  Thus we may argue by induction.  Let $\{\psi_n\}_{n=1}^\infty
  \subset \Mod(S)$ be a family of pseudo-Anosov elements whose
  corresponding closed geodesics on $\calM(S)$ form a dense family in
  $\calM^1(S)$, and let $X \in \Teich(S)$ be a basepoint.  Let
  $\delta_n \to 0$ be given so that the $\delta_n$ neighborhood of the
  axis $A_n$ for $\psi_n$ lies entirely within $\Teich(S)$.  It
  suffices to find a geodesic ray $\ray$ based at $X$ so that for each
  $n$ there is a segment along $\ray$ that lies within $\delta_n$ of the
  axis of some conjugate in $\Mod(S)$ of $\psi_n$ for a full period
  $g_n$ along the axis.

Assume that for $k>1$ we have a ray $\ray_k$ based at $X$ forward
asymptotic to the axis of a conjugate $\hat\psi_k$ of
$\psi_k$ so that the support $\hat \lambda_k^+$ of the attracting
lamination
of  $\hat\psi_k$ lies in the intersection
$$V_k = U_{\delta_1}(\hat{\psi}_1) \cap \ldots 
\cap U_{\delta_{k-1}}(\hat{\psi}_{k-1}).$$
  Then for a sufficiently large power $p_{k+1}$, the support
$\lambda_{k+1}$ of the attracting
  lamination 
for $\psi_{k+1}$ has image
  $\hat{\psi_k}^{p_{k+1}}(\lambda_{k+1})$ within $V_k$.  Taking
  $\ray_{k+1}$ to be the ray asymptotic to the axis of the
  pseudo-Anosov conjugate
$$\hat{\psi}_{k+1} = 
\hat{\psi}_k^{p_{k+1}} \compos \psi_{k+1} \compos
\hat{\psi}_k^{-p_{k+1}}$$ of $\psi_{k+1}$, we have a ray asymptotic to
the axis of a pseudo-Anosov element with attracting lamination in the
intersection $$V_{k+1} = U_{\delta_1}(\hat{\psi}_1) \cap \ldots 
\cap U_{\delta_{k}}(\hat{\psi}_{k}).$$

Thus, $\ray_{k+1}$ lies within $\delta_n$ of
the axis of the conjugate $\hat{\psi}_n$ of $\psi_n$, $n = 1, \ldots,
k+1$, for a full period $g_n$ along the axis of each.  This completes
the induction.

Thus any limit $\ray_\infty$ of $\ray_k$ as $k \to \infty$ in the
visual sphere at $X$ will have a dense trajectory in its projection
$\calM^1(S)$, provided, once again, that it is an infinite ray.  But
$\ray_k$ passes within $\delta_n$ of the axis $\hat A_n$ of
$\hat{\psi}_n$ at the segment $g_n \subset \hat A_n$, for each $k >n$,
so the closest points $\ray_k(t_n)$ to $g_n$ range in a compact
neighborhood in $\Teich(S)$ of a bounded interval along $\hat A_n$.
Thus, given $T>0$, and $n$ so that $t_n > T$, the segments
$\ray_k([0,T])$ sit as subsegments in a family of segments
$\ray_k(t_n)$ whose endpoints converge in $\Teich(S)$ as $k \to
\infty$. Thus, the sequence of geodesics $\ray_k([0,T])$ converges to
a geodesic in Teichm\"uller space for each $T$, by geodesic convexity
of $\Teich(S)$ \cite{Wolpert:Nielsen}.  It follows that the limit
$\ray_\infty$ is infinite and projects to a dense subset of
$\calM^1(S)$ as was claimed.
\end{proof}

\bibliographystyle{math}
\bibliography{math} 

{\sc \small
 \bigskip

\noindent Brown University \bigskip

\noindent University of  Chicago  \bigskip

\noindent Yale University

}

\end{document}